\documentclass[12pt,reqno]{amsart}
\usepackage{amscd}
\usepackage{amssymb}
\usepackage{epsfig}
\usepackage{float}
\usepackage{url}
\usepackage[all]{xy}
\usepackage{graphicx}

\newcommand{\Z}{{\mathbb Z}}
\newcommand{\Q}{{\mathbb Q}}
\newcommand{\C}{{\mathbb C}}

\renewcommand{\P}{{\mathbb P}}

\newcommand{\GG}{{\mathcal G}}

\newcommand{\www}{\widetilde}

\newcommand{\oooo}{\overline}

\DeclareMathOperator{\Aut}{Aut}

\DeclareMathOperator{\End}{End}

\DeclareMathOperator{\id}{id}

\DeclareMathOperator{\lcm}{lcm}

\DeclareMathOperator{\norm}{Norm}
\DeclareMathOperator{\ord}{ord}

\begin{document}

\theoremstyle{plain}
\newtheorem{lemma}{Lemma}[section]
\newtheorem{definition/lemma}[lemma]{Definition/Lemma}
\newtheorem{theorem}[lemma]{Theorem}
\newtheorem{proposition}[lemma]{Proposition}
\newtheorem{corollary}[lemma]{Corollary}
\newtheorem{conjecture}[lemma]{Conjecture}
\newtheorem{conjectures}[lemma]{Conjectures}

\theoremstyle{definition}
\newtheorem{definition}[lemma]{Definition}
\newtheorem{withouttitle}[lemma]{}
\newtheorem{remark}[lemma]{Remark}
\newtheorem{remarks}[lemma]{Remarks}
\newtheorem{example}[lemma]{Example}
\newtheorem{examples}[lemma]{Examples}
\newtheorem{notations}[lemma]{Notations}

\title
[Automorphisms of a $\Z$-lattice with cyclic finite monodromy]
{Automorphisms with eigenvalues in $S^1$ 
of a $\Z$-lattice with cyclic finite monodromy}

\author{Claus Hertling}

\address{Claus Hertling\\Universit\"at Mannheim\\ 
Lehrstuhl f\"ur Mathematik VI\\
Seminargeb\"aude A 5, 6\\
68131 Mannheim, Germany}

\email{hertling@math.uni-mannheim.de}

\thanks{This work was supported by the DFG grant He2287/4-1
(SISYPH)}



\date{January 24, 2018}

\begin{abstract}
For any finite set $M\subset \Z_{\geq 1}$ of positive integers,
there is up to isomorphism a unique $\Z$-lattice $H_M$
with a cyclic automorphism $h_M:H_M\to H_M$ whose eigenvalues
are the unit roots with orders in $M$ and have multiplicity 1.
The paper studies the automorphisms of the pair $(H_M,h_M)$
which have eigenvalues in $S^1$. The main result are 
necessary and sufficient conditions on the set $M$
such that the only such automorphisms are $\pm h_M^k,k\in\Z$. 
The proof uses resultants and cyclotomic polynomials. 
It is elementary, but involved.
Special cases of the main result have been applied to the
study of the automorphisms of Milnor lattices of 
isolated hypersurface singularities.
\end{abstract}

\maketitle

\tableofcontents

\setcounter{section}{0}

\section{Introduction and main result}\label{c1}
\setcounter{equation}{0}

\noindent
In the study of the Milnor lattices of isolated hypersurface singularities,
triples $(H_M,h_M,S)$ with the following properties arise
(e.g. \cite{He11}).

$H_M$ is a $\Z$-lattice of a finite rank $n\in\Z_{\geq 1}$.
It comes with an automorphism $h_M:H_M\to H_M$, which
is called {\it monodromy}, 
and with an $h_M$-invariant bilinear form $S$.
The monodromy is quasiunipotent, 
i.e. its eigenvalues are unit roots,
all eigenvalues have multiplicity 1, 
and $H_M$ has a cyclic generator with respect to $h_M$,
i.e. an element $e_1\in H_M$ with 
\begin{eqnarray}\label{1.1}
H_M = \bigoplus_{i=1}^n \Z\cdot h_M^{i-1}(e_1).
\end{eqnarray}
The restriction of the bilinear form to
the sum $\bigoplus_{\lambda\neq\pm 1}H_\lambda$ is nondegenerate.
Here $H_\C:= H_M\otimes_\Z\C$ and 
$H_\lambda:=\ker(h_M-\lambda\id:H_\C\to H_\C)$ 
is the eigenspace with eigenvalue $\lambda$.

The pair $(H_M,h_M)$ up to isomorphism is determined by the set
\begin{eqnarray}\label{1.2}
M:=\{m\in\Z_{\geq 1}\, |\, e^{2\pi i/m}
\textup{ is an eigenvalue of }h_M\}.
\end{eqnarray}
The characteristic polynomial is $\prod_{m\in M}\Phi_m$.
Here $\Phi_m$ is the cyclotomic polynomial whose zeros are the
unit roots of order $m$.
In the singularity case, the bilinear form may be the 
intersection form or the Seifert form.
For the following problem, the precise form of the
bilinear form $S$ does not matter, only the properties above.
Lemma \ref{t4.1} will show
\begin{eqnarray}\label{1.3}
\Aut(H_M,h_M,S)&=&\{a\in \Aut(H_M,h_M)\, |\, \\
&& \textup{all eigenvalues
of }a\textup{ are in }S^1\}.\nonumber
\end{eqnarray}

The problem is to determine the conditions on the eigenvalues of $h_M$
such that $\Aut(H_M,h_M,S)=\{\pm h_M^k\, |\, k\in\Z\}.$
Theorem \ref{t1.2} gives the complete answer. The conditions are
involved and are given as properties of a graph $\GG(M)$.
The graph and the conditions are formulated in the following definition.

\begin{definition}\label{t1.1}
Let $M\subset \Z_{\geq 1}$ be a finite set of positive integers.

\medskip
(a) A graph $\GG(M)=(M,E(M))$ is associated to it as follows. 
$M$ itself is the set of vertices. The edges in $E(M)$ are directed.
The set $E(m)$ is defined as follows. From a vertex $m_1\in M$ to
a vertex $m_2\in M$ there is no edge if at least one of the following 
two conditions holds:
\begin{list}{}{}
\item[(i)] 
$m_1/m_2$ is not a power of a prime number.
\item[(ii)]
An $m_3\in M-\{m_1,m_2\}$ with $m_2|m_3|m_1$ exists.
\end{list}
If $m_1/m_2$ is a power $p^k$ with $k\in\Z_{\geq 1}$ of a prime number $p$
and if no $m_3\in M-\{m_1,m_2\}$ with $m_2|m_3|m_1$ exists, 
then there is a directed edge from $m_1$ to $m_2$, 
which is additionally labelled with $p$. It is called a $p$-edge. 
Together such edges form the set $E(M)$ of all edges.

\medskip
(b) For any prime number $p$ the components of the graph 
$(M,E(M)-\{p\textup{-edges}\})$ which is obtained by deleting all
$p$-edges, are called the $p$-planes of the graph.
A $p$-plane is called a highest $p$-plane if no $p$-edge ends 
at a vertex
of the $p$-plane. A $p$-edge from $m_1$ to $m_2$ is called a
highest $p$-edge if no $p$-edge ends at $m_1$. 

\medskip
(c) A property $(T_p)$ for a prime number $p$ and a property
$(S_2)$ for the prime number 2:
\begin{eqnarray}\label{1.4}
(T_p) &:& \textup{The graph }\GG(M)
\textup{ has only one highest }p\textup{-plane.}\\
(S_2) &:& \textup{The graph }(M,E(M)-\{\textup{highest }2\textup{-edges}\})
\nonumber \\
&& \textup{ has only 1 or 2 components.}\label{1.5}
\end{eqnarray}

(d) The least common multiple of the numbers
in $M$ is denoted $\lcm(M)\in\Z_{\geq 1}$. 
For any prime number $p$ denote 
\begin{eqnarray*}
l(m,p)&:=&\max(l\in\Z_{\geq 0}\, |\, p^l\textup{ divides }m)
\quad\textup{for any  }m\in\Z_{\geq 1},\\
l(M,p)&:=& \max(l(m,p)\, |\, m\in M)=l(\lcm(M),p).
\end{eqnarray*}
Then $m=\prod_{p\textup{ prime number}} p^{l(m,p)}.$
\end{definition}

The conditions will be discussed after theorem \ref{t1.2}
in the remarks \ref{t1.3}. Examples will be given in \ref{t1.4}.
The following theorem is the main result of this paper.

\begin{theorem}\label{t1.2}
Let $M\subset\Z_{\geq 1}$ be a finite set of positive integers,
and let $(H_M,h_M,S)$ be a triple as above such that 
$M$ is the set of orders of the eigenvalues of $h_M$. Then 
\begin{eqnarray}\label{1.6}
\Aut(H_M,h_M,S)=\{\pm h_M^k\, |\, k\in\Z\}
\end{eqnarray}
holds if and only if the graph $\GG(M)$ satisfies one of the following
two properties.
\begin{list}{}{}
\item[(I)] $\GG(M)$ is connected. It satisfies $(S_2)$.
It satisfies $(T_p)$ for any prime number $p\geq 3$.
\item[(II)] $\GG(M)$ has two components $M_1$ and $M_2$.
The graphs $\GG(M_1)$ and $\GG(M_2)$ are $2$-planes of $\GG(M)$
and satisfy $(T_p)$ for any prime number $p\geq 3$.
Furthermore
\begin{eqnarray}\label{1.7}
\gcd(\lcm(M_1),\lcm(M_2))&\in& \{1;2\},\\
l(M_1,2)&>&l(M_2,2)\in\{0;1\}.\label{1.8}
\end{eqnarray}
\end{list}
\end{theorem}

The theorem will be proved in the sections \ref{c4}, \ref{c5} and 
\ref{c6}. 

\begin{remarks}\label{t1.3}
Let $M$ and $\GG(M)$ be as in definition \ref{t1.1}.

(i) For any $l\in\Z_{\geq 1}$ and any prime number $p$, the set
$\{m\in M\, |\, l(m,p)=l\}$ consists of finitely many 
$p$-planes. 

(ii) From $\GG(M)$ and a prime number $p$, one obtains a
smaller graph $\GG(M)^{(p)}$ as follows. Its vertices are the
$p$-planes of $\GG(M)$. There is a directed edge from a $p$-plane
$E_1$ to a $p$-plane $E_2$ if $E(M)$ contains a $p$-edge
from a vertex in $E_1$ to a vertex in $E_2$.

(iii) The condition $(T_p)$ is equivalent to the condition
that there is a vertex in $\GG(M)^{(p)}$ from which one can reach
all other vertices in $\GG(M)^{(p)}$ if one follows some directed
edges. Especially, $(T_p)$ implies that $\GG(M)$ is connected.

(iv) Any highest $2$-plane is a component of the graph
$(M,E(M)-\{\textup{highest }2\textup{-edges}\})$.
Therefore, if $(S_2)$ holds and $\GG(M)$ is connected,
also $(T_2)$ holds. If $(S_2)$ holds and $\GG(M)$ is not
connected, then $\GG(M)$ has 2 components and each of them
is a 2-plane (and thus there are no 2-edges).
\end{remarks}

\begin{examples}\label{t1.4}
(i) The graph $\GG(M)$ for $M:=\{12,6,4,3,2\}$ has the 
\begin{eqnarray*}
\begin{array}{ll}
\textup{2-edges: } (12,6),(6,3),(4,2),&
\textup{2-planes: } \{12,4\},\ \{6,2\},\ \{3\},\\
\textup{3-edges: } (12,4),(6,2),&
\textup{3-planes: } \{12,6,3\},\ \{4,2\}.
\end{array}
\end{eqnarray*}
For all prime numbers $p\geq 5$, 
$M$ itself is the only $p$-plane, and there are no $p$-edges. 
Case (I) of theorem \ref{t1.2} holds.
The highest 2-edges are $(12,6),(12,3),(4,2)$, and all 3-edges are
highest 3-edges.
The graphs $\GG(M)^{(2)}$ and $\GG(M)^{(3)}$ are just
directed chains with 3 respectively 2 vertices.
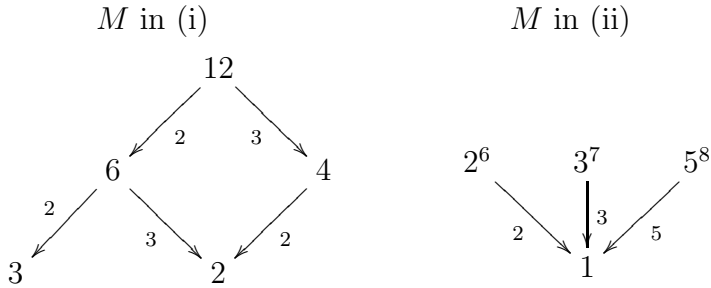
\begin{eqnarray*}
M\textup{ in (i)} \hspace*{4cm}
M\textup{ in (ii)} \\
\begin{xy} \hspace*{-4cm}
\xymatrix{
& & 12 \ar[dl]^2 \ar[dr]_3 &  \\ 
& 6 \ar[dr]_{3} \ar[dl]_{2} & & 4 \ar[dl]^{2} \\ 3 & & 2 &
}\hspace*{6cm}
\xymatrix{ & & \\
2^6 \ar[dr]_{2} & 3^7 \ar[d]^3 & 5^8 \ar[dl]^{5}\\ & 1 &
}\end{xy}
\end{eqnarray*}


(ii) Case (I) of theorem \ref{t1.2} holds also for 
the graph $\GG(M)$ of the set $M=\{2^6,3^7,5^8,1\}$. 
The graph has the 
\begin{eqnarray*}
\begin{array}{ll}
\textup{2-edge: } (2^6,1),&
\textup{2-planes: } \{2^6\},\{3^7,5^8,1\},\\
\textup{3-edge: } (3^7,1),&
\textup{3-planes: } \{3^7\},\{2^6,5^8,1\},\\
\textup{5-edge: } (5^8,1),&
\textup{5-planes: } \{5^8\},\{2^6,3^7,1\}.
\end{array}
\end{eqnarray*}

For all prime numbers $p\geq 7$ $M$ itself is the only $p$-plane,
and there are no $p$-edges. For any $p\in\{2,3,5\}$, the $p$-edge
is a highest $p$-edge.

\medskip
(iii) Case (II) of theorem \ref{t1.2} holds for the graph
$\GG(M)$ of the set $M=\{45,15,14,5,3,2\}$ with 
$M_1=\{14,2\}$ and $M_2=\{45,15,5,3\}$.
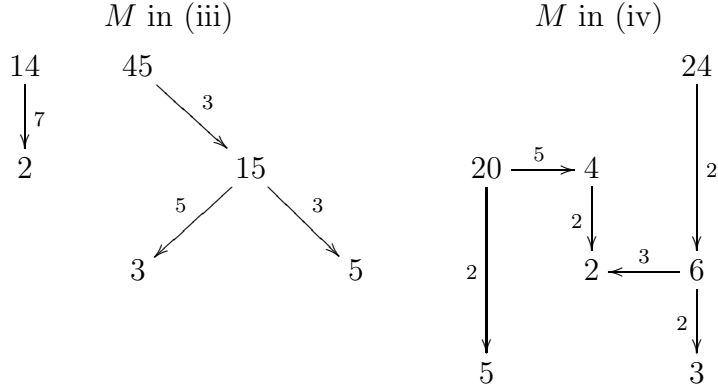
\begin{eqnarray*}
M\textup{ in (iii)} \hspace*{4cm}
M\textup{ in (iv)} \\
\begin{xy} \hspace*{-4cm}
\xymatrix{
14 \ar[d]^7 & 45 \ar[dr]^3 & & \\ 
2 & & 15 \ar[dl]_5 \ar[dr]^3 & \\ 
 & 3 & & 5 
}\hspace*{6cm}
\xymatrix{ & & 24 \ar[dd]^2 \\
20 \ar[dd]_{2} \ar[r]^5 & 4 \ar[d]_2 & \\ 
 & 2 & 6 \ar[l]_3 \ar[d]_2 \\ 5 & & 3 
}\end{xy}
\end{eqnarray*}


(iv) The graph $\GG(M)$ of the set $M=\{24,20,6,5,4,3,2\}$
is connected and satisfies $(T_3)$ and $(T_5)$,
but not $(T_2)$ and thus not $(S_2)$. It has the 
\begin{eqnarray*}
\begin{array}{ll}
\textup{2-planes:}& \{24\},\{20,4\},\{6,2\},\{5\},\{3\}\\
\textup{3-planes:}& \{24,6,3\},\{20,5,4,2\},\\
\textup{5-planes:}& \{20,5\},\{24,6,4,3,2\}.
\end{array}
\end{eqnarray*}
The 2-planes $\{24\}$ and $\{20,4\}$ are both highest 2-planes.
Theorem \ref{t1.2} says 
$\Aut(H_M,h_M,S)\supsetneqq\{\pm h_M^k\, |\, k\in\Z\}$.

\medskip
(v) The graph $\GG(M)$ of the set $M=\{2^6,3^7,5^8\}$
has 3 components. Theorem \ref{t1.2} says 
$\Aut(H_M,h_M,S)\supsetneqq\{\pm h_M^k\, |\, k\in\Z\}$.

\medskip
(vi) Lemma 8.2 in \cite{He11} gives the following
sufficient condition for 
$\Aut(H_M,h_M,S)= \{\pm h_M^k\, |\, k\in\Z\}$.
It is a special case of case (I) in theorem \ref{t1.2}.
$M$ contains a largest number $m_1$ such that $\GG(M)$
is a directed graph with root $m_1$. This implies $(T_p)$
for any $p$. Additionally, a chain of 2-edges exists which
connects all 2-planes. This implies $(S_2)$.
The more special case where $M$ is a 2-plane and 
a directed graph with root $m_1$, was considered
and applied in the proof of \cite[proposition 6.3]{He98}.
\end{examples}

The special case \cite[lemma 8.2]{He11} of theorem \ref{t1.2}
was applied in \cite{He98}, \cite{He11}, \cite{GH16} and \cite{GH17}
in order to study automorphism groups of Milnor lattices of 
isolated hypersurface singularities. Though often not the full Milnor
lattice is a triple $(H_M,h_M,S)$ as above, but it contains
sublattices which are such triples.
In this form, \cite[lemma 8.2]{He11} applies also to some
singularities whose monodromy is not semisimple
(the $T_{pqr}$ in \cite{GH16}). 
But we expect that for other singularities
the more general conditions in theorem \ref{t1.2} will be needed,
and we hope that they will be satisfied.
Conjecture 1.4 in \cite{HZ18} makes our expectations for 
quasihomogeneous singularities precise.

\medskip
The paper is organized as follows. 
Section \ref{c2} studies the resultant of unitary polynomials $f$ and $g$ 
with coefficients in $\Z[x]$ and its relation to the sublattices
$(f,g)\cap\Z[x]_{\leq k}$ of the lattices 
$\Z[x]_{\leq k}:=\{a\in\Z[x]\, |\, \deg a\leq k\}$
for $k\in\Z_{\geq 0}$. Lemma \ref{t2.3} gives fundamental properties,
the lemmata \ref{t2.4} and \ref{t2.5} give statements which will
be applied in the proof of theorem \ref{t1.2}.
Section \ref{c3} recalls in theorem \ref{t3.1} 
basic properties of the cyclotomic polynomials
$\Phi_m$, including the values $\Phi_m(1)$ and Apostol's 
formulas for the resultant of two cyclotomic polynomials \cite{Ap70}.
Theorem \ref{t3.4} gives a tie between different cyclotomic
polynomials which is crucial for the proof of the sufficiency
of the conditions in theorem \ref{t1.2}. It was stated before
as lemma 6.5 in \cite{He98}.
Section \ref{c4} proves the necessity of the conditions in theorem
\ref{t1.2} in the case when $\GG(M)$ is connected.
Section \ref{c5} proves the sufficiency of the conditions in this case.
Section \ref{c6} proves theorem \ref{t1.2} in the case when
$\GG(M)$ is not connected.

\begin{notations}\label{t1.5}
For any polynomial $f\in\C[x]-\{0\}$, the coefficients are denoted
$f_0,...,f_{\deg f}\in\C$. If $f=0$, then $f_0:=0$ and $\deg f:=-\infty$. 

The empty product has value 1. The empty sum has value 0.

$\lambda$ will always denote a unit root in $S^1\subset \C$,
and $\ord(\lambda)$ will be its order, i.e. the minimal $k\in\Z_{\geq 1}$
with $\lambda^k=1$.

$e(z)$ for $z\in\C$ will denote $e^{2\pi i z}\in\C$, 
so for example $e(r)$ for $r\in\Q$ is a unit root.

For $m\in\Z_{\geq 1}$ denote $\Z/m\Z=:\Z_m$, and for
$a\in\Z$ denote its class in $\Z_m$ by $[a]_m$.
\end{notations}

\section{Resultants of unitary polynomials in $\Z[x]$}\label{c2}
\setcounter{equation}{0}

\noindent
The resultant of two polynomials is a very classical object. 
One reference for the following definition is \cite[\S 34]{vW71}.

\begin{definition}\label{t2.1}
The resultant of two polynomials 
$f=\sum_{i=0}^m f_ix^i\in\C[x]-\{0\}$
and $g=\sum_{j=0}^ng_jx^j\in\C[x]-\{0\}$ 
of degrees $\deg f=m,\deg g=n$
with $m+n\geq 1$ is $R(f,g):=\det A(f,g)\in\C$ 
where $A(f,g)\in M((m+n)\times (m+n),\C)$ is
the matrix
\begin{eqnarray}\label{2.1}
A(f,g)=\begin{pmatrix}
f_0 & 0 & \dots & 0 & g_0 & 0 & \dots & 0 \\
f_1 & f_0 & \ddots & \vdots & g_1 & g_0 & \ddots & \vdots \\
\vdots & f_1 & \ddots & 0 & \vdots & g_1 & \ddots & \vdots \\
\vdots & \ddots & \ddots & f_0 & \vdots & \ddots & \ddots & g_0 \\
f_m & \ddots & \ddots & f_1 & g_n & \ddots & \ddots & g_1 \\
0 & \ddots & \ddots & \vdots & 0 & \ddots & \ddots & \vdots \\
\vdots & \ddots & \ddots & \vdots & \vdots & \ddots & \ddots & \vdots \\
0 & \dots & 0 & f_m & 0 & \dots & 0 & g_n 
\end{pmatrix}
\end{eqnarray}
whose first $n$ columns contain the coefficients of $f$
and whose last $m$ columns contain the coefficients of $g$.
In other words, it is the matrix with
\begin{eqnarray}\label{2.2}
(f,xf,...,x^{n-1}f,g,xg,...,x^{m-1}g) = (1,x,...,x^{m+n-1})\cdot A(f,g).
\end{eqnarray}
In the case $m+n=0$ one defines $R(f,g):=1$.
\end{definition}

The basic properties of the resultant are well known.

\begin{proposition}\label{t2.2}
(a) Let $f$ and $g\in\C[x]$ be as in definition \ref{t2.1}.
Let $a_1,...,a_m\in\C$ and $b_1,...,b_n\in\C$ be the zeros of $f$ and $g$, so
$$f=f_0\prod_{i=1}^m(x-a_i),\quad g=g_0\prod_{j=1}^n(x-b_j).$$
Then
\begin{eqnarray}\label{2.3}
R(f,g)&=&f_0^n g_0^m\cdot \prod_{i=1}^m\prod_{j=1}^n(a_i-b_j)\\
&=&(-1)^{m\cdot n}R(g,f),\label{2.4}\\
R(f,g)\neq 0&\iff& \gcd(f,g)_{\C[x]}=1.\label{2.5}
\end{eqnarray}

(b) If $f,g,h\in\C[x]-\{0\}$ then
\begin{eqnarray}\label{2.6}
R(f,gh)=R(f,g)\cdot R(f,h).
\end{eqnarray}
If $f^{(1)},...,f^{(r)},g^{(1)},...,g^{(s)}\in\C[x]-\{0\}$ then
\begin{eqnarray}\label{2.7}
R(\prod_{i=1}^r f^{(i)},\prod_{j=1}^s g^{(j)})
=\prod_{i=1}^r\prod_{j=1}^s R(f^{(i)},g^{(j)}).
\end{eqnarray}
\end{proposition}

\eqref{2.3} is proved for example in \cite[\S 35]{vW71}, \eqref{2.4}, 
\eqref{2.5} and \eqref{2.6} follow from \eqref{2.3}, and \eqref{2.7}
follows from \eqref{2.6} and \eqref{2.4}

We are mainly interested in $R(f,g)$ where $f$ and $g$ are unitary
polynomials (i.e. $f_{\deg f}=1,g_{\deg g}=1$) in $\Z[x]$.
We denote for $k\in\Z_{\geq -1}$
\begin{eqnarray}\label{2.8}
\C[x]_{\leq k}&:=&\{h\in\C[x]\, |, \deg h\leq k\},\\
\Z[x]_{\leq k}&:=&\C[x]_{\leq k}\cap \Z[x]\nonumber
\end{eqnarray}
(so that $\C[x]_{\leq -1}=\Z[x]_{\leq -1}=\{0\}$).

\begin{lemma}\label{t2.3}
Let $f,g\in\Z[x]$ be unitary polynomials of degrees
$m=\deg f,n=\deg g$. They generate
an ideal $(f,g)\subset \Z[x]$ (here $\Z[x]$ is also considered
as an ideal).

(a) 
\begin{eqnarray}\label{2.9}
\Z[x]_{\leq n-1}\cdot f + \Z[x]_{\leq m-1}\cdot g 
=(f,g)\cap \Z[x]_{\leq m+n-1}.
\end{eqnarray}

(b) The $\Z$-lattice in \eqref{2.9} has rank $m+n$ if and
only if $R(f,g)\neq 0$, and then
\begin{eqnarray}\label{2.10}
|R(f,g)| = | \frac{\Z[x]_{\leq m+n-1}}{(f,g)\cap \Z[x]_{\leq m+n-1}}|
\in\Z_{>0}.
\end{eqnarray}

(c) Suppose that $R(f,g)\neq 0$. Then polynomials
$h^{(0)}, h^{(1)},..., h^{(m+n-1)}\in\Z[x]$ with the properties 
in \eqref{2.11}--\eqref{2.14} exist:
\begin{eqnarray}\label{2.11}
(f,g)\cap\Z[x]_{\leq m+n-1}=\bigoplus_{i=0}^{m+n-1}\Z\cdot h^{(i)},\\
\deg h^{(i)}=i\textup{ for any }i,\label{2.12}\\
h_i^{(i)}>0 \textup{ for any }i,\quad
h_{i+1}^{(i+1)}|h_i^{(i)}\textup{ for any }i<m+n-1,\label{2.13}\\
h_i^{(i)}=1 \textup{ for all }i\geq \min(\deg f,\deg g).\label{2.14}
\end{eqnarray}

The coefficients $h_0^{(0)},h_1^{(1)},...,h_{m+n-1}^{(m+n-1)}$ are unique.
\begin{eqnarray}\label{2.15}
|R(f,g)|=\prod_{i=0}^{m+n-1}h_i^{(i)}.
\end{eqnarray}

(d) 
\begin{eqnarray}\label{2.16}
|R(f,g)|=1\iff (f,g)=\Z[x].
\end{eqnarray}
\end{lemma}

{\bf Proof:}
In the case $m=n=0$ $f=g=1$ and all statements are trivial.
So we restrict to the case $m+n\geq 1$.

(a) The ideal is $(f,g)=\Z[x]\cdot f+ \Z[x]\cdot g$.
$\subset$ is obvious. 

$\supset$: For any $h\in (f,g)\cap \Z_{\leq m+n-1}$ 
let $a,b\in\Z[x]$ be such that
$$h=a\cdot f + b\cdot g$$
and such $\max(m+\deg a,n+\deg b)$ is minimal. We have to show
$\max(m+\deg a,n+\deg b)< m+n$.
Suppose that $\max(m+\deg a,n+\deg b)\geq m+n$ for some 
$h\in(f,g)\cap\Z[x]_{\leq m+n-1}$. As $\deg h\leq m+n-1$ and $f$ and
$g$ are unitary, $m+\deg a=n+\deg b$ and 
$a_{\deg a}+b_{\deg b}=0$. Therefore
\begin{eqnarray*}
h=(a-a_{\deg a}\cdot x^{\deg a-n}\cdot g)\cdot f 
+ (b-b_{\deg b}\cdot x^{\deg b-m}\cdot f)\cdot g.
\end{eqnarray*}
Obviously $\deg (a-a_{\deg a}\cdot x^{\deg a-n}\cdot g)<\deg a$ and
$\deg (b-b_{\deg b}\cdot x^{\deg b-m}\cdot f)<\deg b$. This is a contradiction to the
minimality of $\max(m+\deg a,n+\deg b)$. 
Thus $\max(m+\deg a,n+\deg b)<m+n$.

\medskip
(b) This follows immediately from \eqref{2.9}, 
\eqref{2.2} and $R(f,g)=\det A(f,g)$.

\medskip
(c) $R(f,g)\neq 0$ implies $\gcd_{\C[x]}(f,g)=1$. Thus
$(f,g)\cap \Z\supsetneqq\{0\}$ and $(f,g)\cap\Z[x]_{\leq k}$ is a $\Z$-lattice
of rank $k+1$ for any $k\in\Z[x]_{\geq 0}$. 
For $i=0,1,...,m+n-1$, let $h^{(i)}$ be an element of 
$(f,g)\cap \Z[x]_{\leq i}-\Z[x]_{\leq i-1}$
(respectively $(f,g)\cap\Z$ if $i=0$) 
such that $h^{(i)}_i\in\Z$ is positive and minimal.

We show inductively for $i=0,1,...,m+n-1$,
that $h^{(0)},...,h^{(i)}$ is a $\Z$-basis of 
$(f,g)\cap \Z[x]_{\leq i}$.
The case $i=0$ is clear. Suppose it is true for some $i<m+n-1$.
Let $a\in (f,g)\cap\Z[x]_{\leq i+1}$. We have to show
$a\in (f,g)\cap\Z[x]_{\leq i}\oplus\Z\cdot h^{(i+1)}$.
The case $\deg a\leq i$ is trivial, so suppose $\deg a=i+1$.
The minimality of $h_{i+1}^{(i+1)}$ shows $h_{i+1}^{(i+1)}|a_{i+1}$.
Thus 
$$a-\frac{a_{i+1}}{h_{i+1}^{(i+1)}}\cdot h^{(i+1)}
\in (f,g)\cap\Z[x]_{\leq i}.$$
This finishes the inductive proof. The case $i=m+n-1$ is \eqref{2.11}.

\eqref{2.12} and $h_i^{(i)}>0$ hold by definition of $h^{(i)}$. 
Observe $\deg(x\cdot h^{(i)})=i+1$. This and the minimality 
of $h_{i+1}^{(i+1)}$ show $h_{i+1}^{(i+1)}|h_i^{(i)}$.
\eqref{2.14} holds because of $f,g\in(f,g)$ and because they are unitary.
The equations
\begin{eqnarray}\label{2.17}
\prod_{j=0}^i h_j^{(j)} = |\frac{\Z[x]_{\leq i}}{(f,g)\cap\Z[x]_{\leq i}}|
\qquad\textup{for }i=0,...,m+n-1
\end{eqnarray}
hold because $h^{(0)},...,h^{(i)}$ is a $\Z$-basis of 
$(f,g)\cap \Z[x]_{\leq i}$ and because the matrix which 
expresses them
as linear combinations of $x^0,x^1,...,x^i$, is triangular. 
Together they show the uniqueness
of $h_0^{(0)},h_1^{(1)},...,h_{m+n-1}^{(m+n-1)}$.
The case $i=m+n-1$ and \eqref{2.10} give \eqref{2.15}.

\medskip
(d) By \eqref{2.15} and \eqref{2.13}, 
$$|R(f,g)|=1\iff h_0^{(0)}=1
\iff (f,g)=\Z[x]$$
\hfill$\Box$

\bigskip

Lemma \ref{t2.4}, lemma \ref{t2.5} and lemma \ref{t3.2}
will be used in section \ref{c4} in order to show that the
condition $(S_2)$ is necessary in case (I) in theorem \ref{t1.2}.
So, there only the case $p=2$ will be used. Though
lemma \ref{t2.4} is fairly interesting in its own right.

\begin{lemma}\label{t2.4}
Let $f,g\in\Z[x]$ be unitary polynomials of degrees $m=\deg f$, 
$n=\deg g$. Suppose $m\geq n$. 
Let $p$ be a prime number. Consider the following 
four conditions.
\begin{list}{}{}
\item[$(1)$] $|R(f,g)|=p^n$.
\item[$(2)$] $(f,g)=(p,g)$.
\item[$(3)$] $p\in (f,g)$.
\item[$(4)$] $f\in (p,g)$.
\end{list}
Then
\begin{eqnarray}\label{2.18}
(2)\iff (3)\&(4) \iff (1)\&(3)\iff (1)\&(4).
\end{eqnarray}
\end{lemma}

{\bf Proof:}
In the case $m=n=0$ $f=g=1$ and (1)--(4) hold trivially.
So we restrict to the case $m+n\geq 1$.

First we show
\begin{eqnarray}\label{2.19}
\Z[x]_{\leq n-1}\cdot p \oplus
\Z[x]_{\leq m-1}\cdot g=
(p,g)\cap\Z[x]_{\leq m+n-1}.
\end{eqnarray}
$\subset$ is trivial. The proof of $\supset$ is similar to
the one of $\supset$ in \eqref{2.9}:
For any $h\in(p,g)\cap\Z[x]_{\leq m+n-1}$ 
let $a,b\in\Z[x]$ be such that
$$h=a\cdot p+b\cdot g$$ 
and such that $\deg a$ is minimal. We will show
$\deg a\leq n-1$. Suppose $\deg a\geq n$.
Then
$$h=(a-a_{\deg a}\cdot x^{\deg a-n}\cdot g)\cdot p + 
(b+a_{\deg a}\cdot p\cdot x^{\deg a-n})\cdot g,$$
and $\deg (a-a_{\deg a}\cdot x^{\deg a-n}\cdot g)<\deg a$,
a contradiction. Thus $\deg a\leq n-1$.
But now $\deg h\leq m+n-1$  and $\deg g=n$ imply 
immediately $\deg b\leq m-1$. This shows \eqref{2.19}.

\eqref{2.19} implies
\begin{eqnarray}\label{2.20}
|\frac{\Z[x]_{\leq m+n-1}}{(p,g)\cap\Z[x]_{\leq m+n-1}}|=p^n.
\end{eqnarray}

Now the equivalences in \eqref{2.18} will be proved.
$(2)\iff (3)\&(4)$ is trivial. 

If $(1)$ holds, then by \eqref{2.10} and \eqref{2.20}
\begin{eqnarray}\label{2.21}
|\frac{\Z[x]_{\leq m+n-1}}{(f,g)\cap\Z[x]_{\leq m+n-1}}|
=p^n=
|\frac{\Z[x]_{\leq m+n-1}}{(p,g)\cap\Z[x]_{\leq m+n-1}}|.
\end{eqnarray}
Therefore, if $(1)$ holds and if 
one of the two sets $(f,g)\cap\Z[x]_{\leq m+n-1}$
and $(p,g)\cap\Z[x]_{\leq m+n-1}$ is a subset of the other,
they are equal. This shows $(1)\& (3)\iff (1)\&(4)\Rightarrow (2)$.

Finally, if $(2)$ holds, then \eqref{2.10} and \eqref{2.20} 
imply $(1)$.
\hfill $\Box$ 

\bigskip

Lemma \ref{t2.5} will be used in the case $p=2$ in section \ref{c4}
for the treatment of the condition $(S_2)$.

\begin{lemma}\label{t2.5}
Let $p$ be a prime number, and let
$f^{(1)},..., f^{(a)},..., f^{(a+b)}, 
g^{(1)},..., g^{(a)},...,g^{(a+c)}
\in\Z[x]$ (with $a\geq 1,b\geq 0,c\geq 0$) be unitary polynomials.
For $i\in\{1,...,a\}$ define polynomials $h^{(i)}$ by
\begin{eqnarray*}
h^{(i)}:=\left\{\begin{array}{ll}
f^{(i)}& \textup{ if }\deg f^{(i)}< \deg g^{(i)},\\
g^{(i)}& \textup{ if }\deg f^{(i)}\geq \deg g^{(i)}.\end{array}
\right.
\end{eqnarray*}
Suppose 
\begin{eqnarray}\label{2.22}
|R(f^{(i)},g^{(j)})|&=&1\qquad\textup{ for any }i\neq j
\textup{ and for }i=j\geq a+1,\hspace*{1cm}\\
(f^{(i)},g^{(i)})&=&(p,h^{(i)}) \quad\textup{ for }i\in\{1,...,a\}. 
\label{2.23} 
\end{eqnarray}
Then
\begin{eqnarray}\label{2.24}
(\prod_{i=1}^{a+b}f^{(i)},\prod_{j=1}^{a+c}g^{(j)})
=(p,\prod_{i=1}^{a}h^{(i)}).
\end{eqnarray}
\end{lemma}

{\bf Proof:} 
First we consider the special case $b=c=0$ and $h^{(i)}=g^{(i)}$.
Define 
$$f:=\prod_{i=1}^af^{(i)},\qquad g:=\prod_{j=1}^a g^{(j)}.$$
Because of \eqref{2.23}, 
$f^{(i)}$ and $g^{(i)}$ satisfy all conditions in lemma \ref{t2.4},
especially condition $(1)$: 
$$|R(f^{(i)},g^{(i)})|=p^{\deg g^{(i)}}.$$
This and \eqref{2.22} and \eqref{2.7} imply
\begin{eqnarray*}
|R(f,g)|=p^{\deg g},
\end{eqnarray*}
which is condition $(1)$ in lemma \ref{t2.4} for $f$ and $g$.
Because $f^{(i)}$ and $g^{(i)}$ satisfy all conditions in
lemma \ref{t2.4}, and because of \eqref{2.19}, there exist
polynomials $q^{(i)}\in\Z[x]_{\leq m-1}$ and 
$r^{(i)}\in\Z[x]_{\leq \deg g^{(i)}-1}$ with
$$f^{(i)}=q^{(i)}\cdot g^{(i)}+p\cdot r^{(i)}.$$
Therefore 
\begin{eqnarray*}
f &=& \prod_{i=1}^a(q^{(i)}\cdot g^{(i)}+p\cdot r^{(i)}) 
=\left(\prod_{i=1}^aq^{(i)}\right)\cdot g + p\cdot \www r
\end{eqnarray*}
for some polynomial $\www r\in\Z[x]_{\leq \deg g-1}$.
This is condition $(4)$ in lemma \ref{t2.4} for $f$ and $g$.
Therefore $f$ and $g$ satisfy all conditions in lemma \ref{t2.4}.
Condition $(2)$ is \eqref{2.24} in the special case.

Now we consider the general case.
We can suppose that the polynomials are numbered such that
\begin{eqnarray*}
g^{(i)}&=& h^{(i)}\textup{ for }1\leq i\leq d\leq a,\\
f^{(i)}&=& h^{(i)}\textup{ for }d+1\leq i\leq a.
\end{eqnarray*}
Define
\begin{eqnarray*}
a^{(1)}&:=& \prod_{i=1}^d f^{(i)},\quad
a^{(2)}:= \prod_{i=d+1}^a f^{(i)},\quad
a^{(3)}:= \prod_{i=a+1}^b f^{(i)},\quad\\
b^{(1)}&:=& \prod_{i=1}^d g^{(i)},\quad
b^{(2)}:= \prod_{i=d+1}^a g^{(i)},\quad
b^{(3)}:= \prod_{i=a+1}^c g^{(i)},\quad\\
h &:=& \prod_{i=1}^a h^{(i)},\qquad\textup{ thus }
\quad h=b^{(1)}\cdot a^{(2)}.
\end{eqnarray*}
\eqref{2.22}, \eqref{2.7} and \eqref{2.16} tell
$$(a^{(3)},b^{(3)})=\Z[x]\quad\textup{ and }\quad
(b^{(1)},a^{(2)})=\Z[x].$$ 
The special case above tells 
$$(a^{(1)},b^{(1)})=(p,b^{(1)})\quad\textup{ and }
\quad (a^{(2)},b^{(2)})=(p,a^{(2)}).$$
The product of the three ideals $(a^{(i)},b^{(i)})$ for $i=1,2,3$ is
\begin{eqnarray*}
&&(a^{(1)},b^{(1)})\cdot (a^{(2)},b^{(2)}) \cdot (a^{(3)},b^{(3)})\\
&=& (a^{(1)},b^{(1)})\cdot (a^{(2)},b^{(2)})\\
&=& (p,b^{(1)})\cdot (p,a^{(2)}) = (p^2,p\cdot b^{(1)},p\cdot a^{(2)},h)\\
&=& (p,h) \qquad(\textup{ this used }(b^{(1)},a^{(2)})=\Z[x]).
\end{eqnarray*}
The left hand side contains the ideal 
$(a^{(1)}a^{(2)}a^{(3)},b^{(1)}b^{(2)}b^{(3)})$,
thus
\begin{eqnarray}\label{2.25}
(a^{(1)}a^{(2)}a^{(3)},b^{(1)}b^{(2)}b^{(3)})\subset (p,h).
\end{eqnarray}
The special case above also tells
$$|R(a^{(1)},b^{(1)})|=p^{\deg b^{(1)}},\qquad 
|R(a^{(2)},b^{(2)})|=p^{\deg a^{(2)}}.$$
Together with \eqref{2.22} and \eqref{2.7} this implies
\begin{eqnarray}\label{2.26}
|R(a^{(1)}a^{(2)}a^{(3)},b^{(1)}b^{(2)}b^{(3)})|=p^{\deg b^{(1)}}
\cdot p^{\deg a^{(2)}} = p^{\deg h}.
\end{eqnarray}
\eqref{2.19} with $g=h$, $n=\deg g=\deg h$ and $m$ such that
$$m+n=\sum_{i=1}^3(\deg a^{(i)}+\deg b^{(i)})$$ shows
\begin{eqnarray}\label{2.27}
|\frac{\Z[x]_{\leq m+n-1}}{(p,h)\cap\Z[x]_{\leq m+n-1}}|=p^{\deg h}.
\end{eqnarray}
Comparison with \eqref{2.26} and \eqref{2.10} shows
\begin{eqnarray}\label{2.28}
|\frac{\Z[x]_{\leq m+n-1}}{(p,h)\cap\Z[x]_{\leq m+n-1}}|
=|\frac{\Z[x]_{\leq m+n-1}}
{(a^{(1)}a^{(2)}a^{(3)},b^{(1)}b^{(2)}b^{(3)})\cap\Z[x]_{\leq m+n-1}}|
\end{eqnarray}
Together with \eqref{2.25} this gives
$$(p,h)=(a^{(1)}a^{(2)}a^{(3)},b^{(1)}b^{(2)}b^{(3)}),$$
which is \eqref{2.24}.
\hfill $\Box$

\section{Some tie between cyclotomic polynomials}\label{c3}
\setcounter{equation}{0}

\noindent
Recall from the notations \ref{t1.5} that $\lambda$ 
denotes always a unit root in $S^1\subset \C$ 
and that its order is $\ord(\lambda)\in\Z_{\geq 1}$.

For $m\in\Z_{\geq 1}$, the cyclotomic polynomial $\Phi_m$ is the polynomial
\begin{eqnarray}\label{3.1}
\Phi_m(x) := \prod_{\lambda:\ \ord\lambda=m}(x-\lambda),
\end{eqnarray}
whose zeros are the $m$-th primitive unit roots.
It is a unitary and irreducible polynomial in $\Z[x]$ of degree 
$\deg\Phi_m=\varphi(m)\in\Z_{\geq 1}$, where $\varphi:\Z_{\geq 1}\to
\Z_{\geq 1}$ is the Euler phi-function
(see e.g. \cite[Ch 1,2]{Wa82}). Except for the irreducibility, this follows
easily inductively from the formula
\begin{eqnarray}\label{3.2}
x^m-1=\prod_{k|m}\Phi_k.
\end{eqnarray}
Using this formula, one can compute the $\Phi_k$ inductively.
For example for $p$ a prime number and $k,m\in\Z_{\geq 1}$
with $p\not|\, m$
\begin{eqnarray}\label{3.3}
\Phi_{p^{k+1}m}(x)=\Phi_{p^km}(x^p)=\Phi_{pm}(x^{p^{k-1}}) 
\quad\textup{and}\quad
  \Phi_{pm}(x)= \frac{\Phi_m(x^p)}{\Phi_m(x)}.
\end{eqnarray}
Recall (see e.g. \cite[Ch 1,2]{Wa82}) that $\Z[e(\frac{1}{m})]$
is the ring of the algebraic integers within $\Q[e(\frac{1}{m})]$
and that 
\begin{eqnarray}\label{3.4}
\Z[e(\frac{1}{m})]\cap S^1 = \{\pm e(\frac{k}{m})\, |\, k\in\Z\}.
\end{eqnarray}
We will also use the norm
\begin{eqnarray}\label{3.5}
\norm_m:\Z[e(\frac{1}{m})]\to\Z,\quad g(e(\frac{1}{m}))\mapsto
\prod_{\lambda:\ \ord(\lambda)=m}g(\lambda).
\end{eqnarray}
An element of $\Z[e(\frac{1}{m})]$ has norm in $\{\pm 1\}$ if and 
only if it is a unit in $\Z[e(\frac{1}{m})]$.

Part (c) of the following theorem is the main result of \cite{Ap70}.
It gives the resultants of the cyclotomic polynomials.
The proof here is much shorter than that in \cite{Ap70}.

\begin{theorem}\label{t3.1}
(a) $\Phi_m(1)=1$ if $m\geq 2$ and $m$ is not a power of a prime number.
$\Phi_{p^k}(1)=p$ if $p$ is a prime number and $k\in\Z_{\geq 1}$.

(b) $1-\lambda$ is a unit in $\Z[\lambda]$ if and only if $\ord(\lambda)$
is not a power of a prime number and not equal to $1$.

(c) \cite{Ap70} For $m,n\in\Z_{\geq 1}$,
\begin{eqnarray}
R(\Phi_m,\Phi_n)&=& 0\qquad\textup{if }m=n.\label{3.6}\\
R(\Phi_m,\Phi_n)&=& 1\qquad\textup{if neither }\frac{m}{n}\textup{ nor }
\frac{n}{m} \nonumber\\
&& \quad \textup{ is a power of a prime number}.\label{3.7}\\
R(\Phi_{p^kn},\Phi_n)&=& R(\Phi_n,\Phi_{p^kn})=p^{\varphi(n)}
\quad\textup{if }p\textup{ is a prime number}\nonumber\\
&&\quad \textup{ and }k\in\Z_{\geq 1}
\textup{ and }(p,k,n)\neq (2,1,1).\label{3.8}\\
R(1,2)&=&-R(2,1)=2.\label{3.9}
\end{eqnarray}
\end{theorem}

{\bf Proof:}
(a) If $p$ is a prime number and $k\in\Z_{\geq 1}$ then
\begin{eqnarray}\label{3.10}
\Phi_{p^k}(x)&=& x^{(p-1)p^{k-1}}+x^{(p-2)p^{k-1}}+...+x^{p^{k-1}}+1,\\
\textup{so }\Phi_{p^k}(1)&=& p.\label{3.11}
\end{eqnarray}
If one divides both sides of \eqref{3.2} by $\Phi_1=(x-1)$ and 
then puts $x=1$, then one obtains
\begin{eqnarray}\label{3.12}
m=\prod_{k:\, k|m,k\neq 1}\Phi_k(1).
\end{eqnarray}
This and \eqref{3.11} and induction show $\Phi_m(1)=1$ for any 
$m\in\Z_{\geq 2}$ which is not a power of a prime number.

\medskip
(b) Let $\lambda$ be a unit root with order $\ord(\lambda)=m$.
\begin{eqnarray}\label{3.13}
\norm_m(1-\lambda)=\prod_{\mu:\, \ord(\mu)=m}(1-\mu)=\Phi_m(1).
\end{eqnarray}
This and part (a) show that $1-\lambda$ is a unit in $\Z[\lambda]$
if and only if $m$ is not a power of a prime number
and not equal to 1.

\medskip
(c) \eqref{2.3} for any $m,n\in\Z_{\geq 1}$ gives
\begin{eqnarray}
R(\Phi_m,\Phi_n) &=& \prod_{\lambda:\, \ord(\lambda)=m}
\prod_{\mu:\, \ord(\mu)=n}(\lambda-\mu)\nonumber \\
&=& \prod_{\lambda:\, \ord(\lambda)=m}
\prod_{\mu:\, \ord(\mu)=n}
\left[\lambda\cdot(1-\lambda^{-1}\cdot\mu)\right]. \label{3.14}
\end{eqnarray}
\eqref{3.6} and \eqref{3.9} follow immediately.
For $m$ and $n$ as in \eqref{3.7}, $\lambda^{-1}\mu$ is a unit root
whose order is not a power of a prime number. Then by (b)
all factors in the product above are units in $\Z[\lambda,\mu]$, 
so the product is a unit in $\Z$, so it is in $\{\pm 1\}$. 
As $(1-\lambda^{-1}\cdot\mu)\cdot (1-\oooo\lambda^{-1}\cdot\oooo\mu)>0$,
the product is positive, thus it is $+1$. This shows \eqref{3.7}.

For $m=p^kn$ with $p,k,n$ as in \eqref{3.8}, write
$\lambda^{-1}=e(\frac{a}{p^kn})$ and $\mu=e(\frac{b}{n})$
with $a\in\{1,...,p^kn\},b\in\{1,...,n\}$ 
with $\gcd(a,p^kn)=1,\gcd(b,n)=1$, so that
$$\lambda^{-1}\mu=e(\frac{a}{p^kn}+\frac{b}{n})=e(\frac{a+p^kb}{p^kn}).$$
Write $n=p^l\cdot c$ with $\gcd(p,c)=1,l\in\Z_{\geq 0}$.
As $\gcd(a+p^kb,p)=\gcd(a,p)=1$, 
the order of $\lambda^{-1}\mu$ is $p^{k+l}\cdot c/\gcd(a+p^kb,c)$.
It is a power of a prime number (namely $p^{k+l}$) 
if and only if $[a+p^kb]_c=0$. How often does this hold?

If $a$ runs through 
$\{\www a\, |\, 1\leq a\leq p^kn,\gcd(\www a,p^kn)=1\}$, 
then $[a]_c$ runs with multiplicity $\varphi(p^kn)/\varphi(c)$ 
through all units in $\Z_c$. 
If $b$ runs through $\{\www b\, |\, 1\leq b\leq n,\gcd(b,n)=1\}$,
then $[b]_c$ and $[p^kb]_c$ run with multiplicity $\varphi(n)/\varphi(c)$
through all units in $\Z_c$.
Thus the sum $[a]_c+[p^kb]_c=[a+p^kb]_c$ vanishes in 
$$\frac{\varphi(p^kn)}{\varphi(c)}\cdot\frac{\varphi(n)}{\varphi(c)}
\cdot \varphi(c) =\varphi(p^{k+l})\cdot \varphi(n)$$
cases. Therefore the product in \eqref{3.14} contains
$\varphi(p^{k+l})\cdot \varphi(n)$ factors $(1-\lambda^{-1}\mu)$
with $\ord(\lambda^{-1}\mu)$ a power of a prime number,
and this power is $p^{k+l}$. Together with \eqref{3.13}
for $p^{k+l}$ instead of $m$ and with
$\Phi_{p^{k+l}}(1)=p$ (part (a)) this shows that these factors
give $p^{\varphi(n)}$. The other factors together give $\pm 1$.
The same argument as above with the complex conjugate unit roots 
shows $R(\Phi_{p^kn},\Phi_n)>0$ if $(p,k,n)\neq (2,1,1)$. 
This proves \eqref{3.8}.
\hfill $\Box$

\bigskip

Lemma \ref{t3.2} will be used in the case $p=2$ in section \ref{c4}
for the treatment of the condition $(S_2)$.

\begin{lemma}\label{t3.2}
Let $p$ be a prime number, let $m\in\Z_{\geq 1}$ with
$p\not|\, m$, and let $k,l_1,...,l_r\in\Z_{\geq 0}$
with $k>l_1>...>l_r$ for some $r\in\Z_{\geq 1}$.
Define 
\begin{eqnarray*}
f:= \Phi_{p^km},\quad 
g:=\Phi_{p^{l_1}m}\cdot ...\cdot\Phi_{p^{l_r}m}.
\end{eqnarray*}
Then $f$ and $g$ satisfy all properties (1)--(4) in lemma
\ref{t2.4}, especially $(f,g)=(p,g)$.
\end{lemma}

{\bf Proof:}
Because of lemma \ref{t2.4}, it will be sufficient to show
$\deg f\geq \deg g$ and  the properties 
(1) and (3) in lemma \ref{t2.4}.

$\deg f\geq \deg g$: 
\begin{eqnarray*}
\deg f&=& \varphi(p^km)=\varphi(p^k)\cdot \varphi(m)
=(p-1)p^{k-1}\cdot \varphi(m),\\
\deg g&\leq & \deg(\Phi_{p^{k-1}m}\cdot \Phi_{p^{k-2}m}
\cdot ...\cdot \Phi_{p^0m})\\
&=& ((p-1)(p^{k-2}+p^{k-3}+...+1)+1)\cdot\varphi(m)\\
&=& p^{k-1}\cdot \varphi(m)\leq \deg f.
\end{eqnarray*}

Property (1), $|R(f,g)|=p^{\deg g}$: This uses \eqref{3.8}
and possibly \eqref{3.9} (if $p=2$ and $m=1$).
\begin{eqnarray*}
|R(f,g)| = \prod_{i=1}^r |R(\Phi_{p^km},\Phi_{p^{l_i}m}|
= \prod_{i=1}^r p^{\varphi(p^{l_i}m)} = p^{\deg g}.
\end{eqnarray*}

Property (3), $p\in (f,g)$: 
$f$ divides 
$$\www f:=\prod_{a|m}\Phi_{p^ka}=\frac{x^{p^km}-1}{x^{p^{k-1}m}-1}
= x^{(p-1)p^{k-1}m}+x^{(p-2)p^{k-1}m}+...+1,$$
and $g$ divides
$$\www g:=\prod_{b|p^{k-1}m}\Phi_{b}=x^{p^{k-1}m}-1,$$
thus $(f,g)\supset (\www f,\www g)$.
Observe
$$ \www f\equiv p\mod (\www g),$$
so $(p,\www g)= (\www f,\www g)\subset (f,g)$. This shows (3).
\hfill$\Box$

\bigskip

Lemma \ref{t3.3} will be used in the proof of theorem \ref{t3.4}.

\begin{lemma}\label{t3.3}
(a) Let $p$ be a prime number and $k,m\in\Z_{\geq 1}$. Then
\begin{eqnarray}\label{3.15}
\Phi_{p^km}(e({1\over m}))=p \cdot \textup{unit}.
\end{eqnarray}
Here and in the proof {\rm unit} means an invertible element
in $\Z[\lambda]$ for a suitable unit root $\lambda$.

\medskip
(b) Let $\lambda$ be a unit root and $m=\ord(\lambda)$ its order. The set
$\{\norm_m(1-\lambda^k)\ |\ k\in \Z\}$ is the union of the set $\{0\}$, the set
$$\{p^{\frac{\varphi (m)}{\varphi(p^l)}} \ |\ l\geq 1\textup{ and } p
\mbox{ a prime number such that }p^l|m\} ,$$
and, if and only if $m$ is not a power of a prime number,
the set $\{1\}$. 
\end{lemma}

{\bf Proof:} (a) 
If $p_i$ are different prime numbers and $k_i\geq 1$, then
\begin{eqnarray}\label{3.16}
\Phi_{p_1^{k_1}...p_l^{k_l}} (x)
= \Phi_{p_1...p_l}(x^{p_1^{k_1-1}...p_l^{k_l-1}}).
\end{eqnarray}
Hence \eqref{3.15} can be reduced to the statement
\begin{eqnarray}\label{3.17}
\Phi_{p_1...p_l}(e({1\over p_2...p_l}))= p_1 \cdot \mbox{ unit}.
\end{eqnarray}
If $p,q$ are prime numbers and if they and $m\in\Z_{\geq 1}$ 
are such that $p\neq q$ and $p$ and $q$ do not divide $m$, then
$$\Phi_{pm} (e({1\over qm})) 
= \prod_{\ord(\lambda)=pm}(e(\frac{1}{qm}-\lambda)
= e({\varphi (pm)\over qm})\cdot 
  \prod_{\ord (\lambda )=pm} (1-\lambda \cdot e({-1\over qm}))$$
is a unit, because the order $\ord(\lambda \cdot e({-1\over qm}))$
is not a power of a prime number.
Using $\Phi_{pm}(x^q)=\Phi_{pqm}(x)\Phi_{pm}(x)$, we get
$$\Phi_{pqm}(e({1\over qm}))=\Phi_{pm}(e({1\over m}))\cdot \mbox{ unit}.$$
Thus \eqref{3.17} can be reduced to the trivial case $\Phi_p(1)=p$.

\medskip
(b) If $\ord (\lambda^k)$ is not a power of a prime number, then
$\norm_m(1-\lambda^k)=1$ because $\Phi_{\ord (\lambda^k)}(1)=1$.
If $\ord (\lambda^k)=p^l $ then
$$\norm_m(1-\lambda^k)=(\Phi_{p^l}(1))^{\frac{\varphi(m)}{\varphi(p^l)}}
=p^{\frac{\varphi(m)}{\varphi(p^l)}}.$$
\hfill $\Box$

\bigskip

Theorem \ref{t3.4} gives a tie between different cyclotomic
polynomials. It will be crucial for the proof in section \ref{c5}
of the sufficiency of the conditions in case (I) in 
theorem \ref{t1.2}. It was stated before as lemma 6.5
in \cite{He98}.

\begin{theorem}\label{t3.4}
Let $p$ be a prime number, $k,m\in \Z_{\geq 1},\ c(x)\in \Z[x]$ such that
$c(e({1\over p^km}))=1$ and $|c(e({1\over m}))|=1$.
\begin{list}{}{}
\item[(a)] If $p\geq 3$ then $c(e({1\over m}))=1$.
\item[(b)] If $p=2$ then $c(e({1\over m}))=\pm 1$.     
\item[(c)] If $p=2$ and $c(e({1\over p^lm}))=1$ for some $l\neq k$ then     
          $c(e({1\over m}))=1$.     
\end{list}
\end{theorem}

{\bf Proof:}
(a) Let $p,k,m,c(x)$ be as in the theorem, with $p\geq 3$.
There exists a polynomial $r(x)\in \Z[x]$ such that
$1-c(x)=\Phi_{p^km}(x)\cdot r(x).$ Then
\begin{eqnarray*}
\norm_m(1-c(e({1\over m})))&=&
  \norm_m(\Phi_{p^km}(e({1\over m})))\cdot
  \norm_m(r(e({1\over m})))\\
&=& (\pm 1)\cdot p^{\varphi (m)}  \cdot
  \norm_m(r(e({1\over m}))).
\end{eqnarray*}
The second equality uses lemma \ref{t3.3} (a).
From \eqref{3.4} and $|c(e({1\over m}))|=1$ we obtain 
$$c(e({1\over m})) \in \{\pm e({l\over m})\ |\ l\in \Z\}.$$

\medskip
{\bf Case 1,} $m$ is odd: Then $\Z[e({1\over m})]=\Z[e({1\over 2m})]$, 
$\norm_m=\norm_{2m}$, \\
$\{\pm e({l\over m})\ |\ l\in \Z\}=\{e({l\over 2m})\ |\ l\in \Z\}$.
Because of lemma \ref{t3.3} (b) and $\varphi(p^l)>1$ for $l\geq 1$, 
the only number in $\{\norm_m(1-e(\frac{l}{2m}))\ |\ l\in \Z\}$, 
which is divisible by $p^{\varphi (2m)}=p^{\varphi (m)}$, is 0. Thus
$\norm_m(1-c(e(\frac{1}{m})))=0$ and $c(e({1\over m}))=1$.

\medskip
{\bf Case 2,} $m$ is even: Then 
$\{\pm e({l\over m})\ |\ l\in \Z\}=\{e({l\over m})\ |\ l\in \Z\}$.
Because of lemma \ref{t3.3} (b) and $\varphi(p^l)>1$ for $l\geq 1$, 
the only number in $\{\norm_m(1-e(\frac{l}{m}))\ |\ l\in \Z\}$,
which is divisible by $p^{\varphi (m)}$, is 0. Thus
$\norm_m(1-c(e(\frac{1}{m})))=0$ and $c(e({1\over m}))=1$.

This proves part (a).

\medskip
(b) Let $p,k,m,c(x)$ be as in the theorem, with $p=2$.
The proof proceeds as the proof of part (a). Only the statement
$\varphi(p^l)>1$ becomes wrong if $l=1$. Then
$\norm_m(1-e(\frac{l}{2m}))=2^{\varphi(m)}$ in case 1
respectively 
$\norm_m(1-e(\frac{l}{m}))=2^{\varphi(m)}$ in case 2
is possible, but only in the case $e(\frac{l}{2m})=-1$ in case 1
respectively $e(\frac{l}{m})=-1$ in case 2,
as the proof of lemma \ref{t3.3} (b) shows.

\medskip
(c) Let $p,k,m,c(x)$ be as in the theorem, with $p=2$
and $c(e(\frac{1}{p^lm}))=1$ for some $l\neq k$.
There exists a polynomial $r(x)\in\Z[x]$ such that
$$1-c(x)=\Phi_{p^km}(x)\cdot \Phi_{p^lm}(x)\cdot r(x).$$
Then
\begin{eqnarray*}
&&\norm_m(1-c(e({1\over m})))\\
&=& \norm_m(\Phi_{p^km}(e({1\over m})))
\cdot \norm_m(\Phi_{p^lm}(e({1\over m})))
\cdot \norm_m(r(e({1\over m})))\\
  &=& (\pm 1)\cdot p^{\varphi(m)}\cdot (\pm 1)\cdot p^{\varphi(m)}
\cdot \norm_m(r(e({1\over m}))).
\end{eqnarray*}
The last equality uses lemma \ref{t3.3} (a).
Now one has again to go through the two cases 
and apply lemma \ref{t3.3} (b).
As $2\varphi(m)$ is bigger than $\varphi(m)/\varphi(p^l)$ in any case,
$\norm_m(1-c(e({1\over m})))=0$ and $c(e({1\over m}))=1$.
\hfill $\Box$

\section{Necessity of the conditions in the main result
in the connected case}\label{c4}
\setcounter{equation}{0}

\noindent
Let $(H_M,h_M,S)$ be a triple as in the introduction
and let $M$ be the set of orders of the eigenvalues of $h_M$.
The main point in this section is the proof that the conditions in case
(I) in theorem \ref{t1.2} are necessary for 
$\Aut(H_M,h_M,S)=\{\pm h_M^k\, |\, k\in\Z\}$ 
if $\GG(M)$ is connected. 
But before, the next lemma shows that the precise form 
of the bilinear form $S$ in the triple $(H_M,h_M,S)$ is unimportant.
Recall that $\prod_{m\in M}\Phi_m$ is the characteristic polynomial
of $h_M$ and that $(\prod_{m\in M}\Phi_m)$ denotes its ideal in $\Z[x]$.

\begin{lemma}\label{t4.1}
\begin{eqnarray}\label{4.1}
\End(H_M,h_M)&=&\{c(h_M)\, |\, c(x)\in\Z[x]\},\\
\Aut(H_M,h_M,S)&=&\{c(h_M)\, |\, c(x)\in\Z[x]\textup{ with }|c(\lambda)|=1
\label{4.2}\\
&& \hspace*{1cm}\textup{ for any eigenvalue }
\lambda\textup{ of }h_M\},\nonumber
\end{eqnarray}
Thus 
\begin{eqnarray}
&&\Aut(H_M,h_M,S)=\{\pm h_M^k\, |\, k\in\Z\}\nonumber\\
&\iff& \{c(x)\in\Z[x]\, |\, |c(\lambda)|=1
\textup{ for any eigenvalue }\lambda\textup{ of }M\} \nonumber\\
&&=\{\pm x^k\, |\, k\in\Z\}+(\prod_{m\in M}\Phi_m).\label{4.3}
\end{eqnarray}
\end{lemma}

{\bf Proof:} Due to \eqref{1.1}, 
for any $B\in\End(H_M,h_M)$ a unique polynomial
$b(x)=\sum_{i=0}^{n-1}b_ix^i\in\Z[x]$ with $B(e_1)=\sum_{i=0}^{n-1}
b_i h_M^i(e_1)$ exists. The commutativity $h_M\circ B=B\circ h_M$ implies
$B=b(h_M)$. This proves \eqref{4.1}. 

Any eigenspace $H_\lambda$ of $h_M$ is 1-dimensional by hypothesis. 
Two eigenspaces $H_\lambda$ and $H_\mu$ are orthogonal with respect to $S$ if 
$\mu\neq\oooo{\lambda}$, because $S$ is $h_M$-invariant.
By hypothesis, the restriction of $S:H_\lambda\times H_{\oooo\lambda}\to\C$
is nondegenerate if $\lambda\notin\{\pm 1\}$.

Now consider an automorphism $b(h_M)\in\Aut(H_M,h_M,S)$ for some 
$b(x)\in\Z[x]$. The space $H_\lambda$ is also an eigenspace of $b(h_M)$,
and it has eigenvalue $b(\lambda)$ on $H_\lambda$. 
As $b(h_M)$ is an automorphism of $H_M$, 
its eigenvalue on $H_1$ if $H_1\neq\{0\}$ 
and its eigenvalue on $H_{-1}$ if $H_{-1}\neq \{0\}$
must be in $\{\pm 1\}$.
It respects $S$ on $\bigoplus_{\lambda\neq\pm 1}H_\lambda$ 
if and only if $|b(\lambda)|=1$ for any eigenvalue $\lambda\neq\pm 1$.
Therefore $|b(\lambda)|=1$ for any eigenvalue.

Vice versa, suppose that $b(h_M)\in\End(H_M,h_M)$ for some $b(x)\in\Z[x]$
with $|b(\lambda)|=1$ for any eigenvalue $\lambda$ of $h_M$.
Then $b(h_M)$ respects $S$, and $\det b(h_M)\in\{\pm 1\}$, and thus
$b(h_M)\in\Aut(H_M,h_M,S)$. This completes the proof of \eqref{4.2}.

\eqref{4.3} is an immediate consequence of \eqref{4.2}
\hfill $\Box$

\bigskip

Suppose now that $\GG(M)$ is connected. We will show
$\supsetneqq$ in \eqref{4.3} if $(T_p)$ does not hold for some
prime number $p\geq 3$ (1st case) or if $(S_2)$ does not hold (2nd case).

\medskip
{\bf 1st case,} $(T_p)$ does not hold for some prime number $p\geq 3$:
Let $E_1,...,E_r$ with $r\geq 2$ be the highest $p$-planes.
Let $F_1\subset M$ be the union of all $p$-planes which can
be reached within the graph $\GG(M)^{(p)}$ 
(whose vertices are all the $p$-planes in $\GG(M)$, see remark \ref{t1.3}(ii))
by starting at $E_1$ and following some directed edges.
Let $F_2$ be the union of all $p$-planes which can
be reached within the graph $\GG(M)^{(p)}$ 
by starting at one of the points $E_2,...,E_r$ 
and following some directed edges.
As $E_1,...,E_r$ are all highest $p$-planes, $F_1\cup F_2=M$.
As $\GG(M)$ is connected, $F_1\cap F_2\neq\emptyset$.
Define 
$$G_1:=F_1-F_1\cap F_2,\ G_2:=F_2-F_1\cap F_2,\ G_3:= F_1\cap F_2,$$
so that $G_1 \dot\cup G_2 \dot\cup G_3=M$.
Also 
$$(F_1\cap F_2)\cap E_1=\emptyset\quad\textup{and}\quad
(F_1\cap F_2)\cap (E_2\cup...\cup E_r)=\emptyset$$
are obvious, and they imply
\begin{eqnarray}\label{4.4}
l(G_3,p)&<& l(E_1,p)=l(F_1,p)=l(G_1,p),\\
l(G_3,p)&<& l(E_2\cup ...\cup E_r,p)
=l(F_2,p)=l(G_2,p).\nonumber
\end{eqnarray}

By definition of $G_1$ and $G_2$, there are no edges at all between
vertices in $G_1$ and vertices in $G_2$. With \eqref{3.7} and \eqref{2.7},
the resultant of the following polynomials is in $\{\pm 1\}$,
\begin{eqnarray*}
R(\prod_{m\in G_1}\Phi_m,\prod_{m\in G_2}\Phi_m)=\pm 1.
\end{eqnarray*}
By \eqref{2.16}, there exist $a_1,a_2\in\Z[x]$ with
\begin{eqnarray*}
1&=& a_1\cdot\prod_{m\in G_1}\Phi_m + a_2\cdot\prod_{m\in G_2}\Phi_m.
\end{eqnarray*}
Write $d:=\lcm(G_3)$ (so that $l(G_3,p)=l(d,p)$). Then
\begin{eqnarray*}
x^d-1&=& a_1\cdot (x^d-1) \cdot\prod_{m\in G_1}\Phi_m
+ a_2\cdot(x^d-1)\cdot\prod_{m\in G_2}\Phi_m\\
&=& b_1\cdot\prod_{m\in G_1\cup G_3}\Phi_m
+ b_2\cdot\prod_{m\in G_2\cup G_3}\Phi_m
\end{eqnarray*}
for some $b_1,b_2\in\Z[x]$. Define 
\begin{eqnarray}\label{4.5}
c(x):= x^d-b_1\cdot\prod_{m\in G_1\cup G_3}\Phi_m
=1+ b_2\cdot\prod_{m\in G_2\cup G_3}\Phi_m\in\Z[x].
\end{eqnarray}
We want to show that there do not exist $\varepsilon\in\{\pm 1\}$
and $k\in\Z_{\geq 0}$ with 
\begin{eqnarray}\label{4.6}
c(x)\equiv \varepsilon\cdot x^k\mod (\prod_{m\in M}\Phi_m).
\end{eqnarray}
We suppose that 
$\varepsilon\in\{\pm 1\}$ and $k\in\Z$ with \eqref{4.6} exist. 
We want to arrive at a contradiction.
\eqref{4.5} and \eqref{4.6} give
\begin{eqnarray*}
e(\frac{d}{m})&=&c(e(\frac{1}{m}))=\varepsilon \cdot e(\frac{k}{m})
\qquad\textup{for }m\in G_1\cup G_3,\\
1&=&c(e(\frac{1}{m}))=\varepsilon \cdot e(\frac{k}{m})
\qquad\textup{for }m\in G_2\cup G_3,
\end{eqnarray*}
thus
\begin{eqnarray*}
d&\equiv& \frac{m}{2}\cdot \delta_{-1,\varepsilon}+k\mod m\Z
\qquad\textup{for } m\in G_1\cup G_3,\\
0&\equiv& \frac{m}{2}\cdot \delta_{-1,\varepsilon}+k\mod m\Z
\qquad\textup{for } m\in G_2\cup G_3.
\end{eqnarray*}
If $\varepsilon=-1$ this shows that $m\in G_1\cup G_2\cup G_3$ is even.
Recall $m=\prod_{q\textup{ prime number}}q^{l(m,q)}$. 
In any case, whether $\varepsilon=-1$ or $\varepsilon=1$, 
\begin{eqnarray*}
d&\equiv& k\mod p^{l(m,p)}\Z
\qquad\textup{for } m\in G_1\cup G_3,\\
0&\equiv& k\mod p^{l(m,p)}\Z
\qquad\textup{for } m\in G_2\cup G_3,\\
\end{eqnarray*}
as $p\geq 3$, so
\begin{eqnarray*}
d&\equiv& k\mod p^{l(G_1,p)}\Z\\
0&\equiv& k\mod p^{l(G_2,p)}\Z
\end{eqnarray*}
As $l(G_1,p)>l(G_3,p)=l(d,p)$ and $l(G_2,p)>l(G_3,p)$ by \eqref{4.4},
\begin{eqnarray*}
d\equiv k\equiv 0\mod p^{l(d,p)+1}.
\end{eqnarray*}
But this is impossible, as it contradicts the definition of $l(d,p)$.
Therefore $\varepsilon\in\{\pm 1\}$ and $k\in\Z$ with \eqref{4.6}
do not exist. Thus $c(h_M)\notin\{\pm h_M^k\, |\, k\in\Z\}$.
On the other hand, \eqref{4.5} and lemma \ref{t4.1} tell
$c(h_M)\in \Aut(H_M,h_M,S)$. This proves the necessity of
$(T_p)$ for $p\geq 3$ in theorem \ref{t1.2} in the case when
$\GG(M)$ is connected.

\medskip
{\bf 2nd case,} $(S_2)$ does not hold:
Let $E_1,...,E_r$ with $r\geq 3$ be the components of the graph
$(M,E(M)-\{\textup{highest }2\textup{-edges}\})$.
As $\GG(M)$ has at least one highest 2-plane, by remark \ref{t1.3} (iv)
we can suppose that $E_1$ is a highest 2-plane and that
$l(E_1,2)=l(M,2)$.
As $\GG(M)$ is connected, we can also suppose that $E_2$ is a component
such that there exists a highest 2-edge from a vertex in $E_1$ to a
vertex in $E_2$. Define $F_1:=E_1\cup E_2\neq\emptyset$ and 
$F_2:=E_3\cup...\cup E_r\neq\emptyset$.
Then $M=F_1\dot\cup F_2$. Consider for any odd $a\in\Z_{\geq 1}$ the sets
\begin{eqnarray*}
B_a&:=&\{2^ka\, |\, k\in\Z_{\geq 0}\},\quad 
B_{1,a}:=B_a\cap F_1,\quad B_{2,a}:=B_a\cap F_2.
\end{eqnarray*}
The construction of $E_1,...,E_r$ and of $F_1$ and $F_2$ tells
\begin{eqnarray}
&&B_{1,a}\neq\emptyset\quad\textup{and}\quad B_{2,a}\neq\emptyset \nonumber\\
&\Rightarrow& \left\{\begin{array}{llll}
\textup{either}&B_{1,a}=\{m_1\}&\textup{ and }&m_1>m\ \forall\ m\in B_{2,a}\\
\textup{or}&B_{2,a}=\{m_1\}&\textup{ and }&m_1>m\ \forall\ m\in B_{1,a}.
\end{array}\right. \label{4.7}
\end{eqnarray}
Define
\begin{eqnarray*}
A_{12}&:=& \{a\in\Z_{\geq 1}\, |\, a\textup{ odd, }B_{1,a}\neq \emptyset,
B_{2,a}\neq\emptyset\},\\
A_1&:=& \{a\in\Z_{\geq 1}\, |\, a\textup{ odd, }B_{1,a}\neq \emptyset,
B_{2,a}=\emptyset\},\\
A_2&:=& \{a\in\Z_{\geq 1}\, |\, a\textup{ odd, }B_{1,a}= \emptyset,
B_{2,a}\neq\emptyset\}.
\end{eqnarray*}
Define 
\begin{eqnarray*}
f_a&:=&\prod_{m\in B_{1,a}}\Phi_m\quad\textup{for }a\in A_{12}\cup A_1,\\
g_a&:=&\prod_{m\in B_{2,a}}\Phi_m\quad\textup{for }a\in A_{12}\cup A_2.
\end{eqnarray*}
If $a\in A_{12}$, then either the pair $(f_a,g_a)$ or the pair
$(g_a,f_a)$ satisfies the properties of the pair $(f,g)$ in lemma \ref{t3.2}
with $p=2$, because of \eqref{4.7}.
Furthermore, observe that the sets of vertices $F_1$ and $F_2$ 
are connected only by some highest 2-edges and not by any other edges.
This implies $|R(f_{a_1},g_{a_2})|=1$ for $a_1\neq a_2$ by \eqref{3.7}. 
Therefore the polynomials $f_a,a\in A_{12}\cup A_1$, and the polynomials
$g_a,a\in A_{12}\cup A_2$, satisfy all properties of the polynomials
$f^{(1)},...,f^{(a+b)},g^{(1)},...,g^{(a+c)}$ in lemma \ref{t2.5}, 
with the obvious differences in the notations.
\eqref{2.24} in lemma \ref{t2.5} tells that there exist polynomials
$b_1,b_2\in\Z[x]$ with
\begin{eqnarray*}
2&=& b_1\cdot \prod_{a\in A_{12}\cup A_1}f_a 
+ b_2\cdot \prod_{a\in A_{12}\cup A_2}g_a\\
&=& b_1\cdot \prod_{m\in F_1}\Phi_m
+ b_2\cdot \prod_{m\in F_2}\Phi_m.
\end{eqnarray*}
Now define
\begin{eqnarray}\label{4.8}
c(x)&:=& 1-b_1\cdot \prod_{m\in F_1}\Phi_m
= -1+b_2\cdot \prod_{m\in F_2}\Phi_m.
\end{eqnarray}
The rest of the argument is similar to the 1st case.
We want to show that there do not exist $\varepsilon\in\{\pm 1\}$
and $k\in\Z_{\geq 0}$ with 
\begin{eqnarray}\label{4.9}
c(x)\equiv \varepsilon\cdot x^k\mod (\prod_{m\in M}\Phi_m).
\end{eqnarray}
We suppose that 
$\varepsilon\in\{\pm 1\}$ and $k\in\Z$ with \eqref{4.9} exist. 
We want to arrive at a contradiction.
\eqref{4.8} and \eqref{4.9} give
\begin{eqnarray*}
1&=&c(e(\frac{1}{m}))=\varepsilon \cdot e(\frac{k}{m})
\qquad\textup{for }m\in F_1,\\
-1&=&c(e(\frac{1}{m}))=\varepsilon \cdot e(\frac{k}{m})
\qquad\textup{for }m\in F_2,
\end{eqnarray*}
thus
\begin{eqnarray*}
0&\equiv& \frac{m}{2}\cdot \delta_{-1,\varepsilon}+k\mod m\Z
\qquad\textup{for } m\in F_1,\\
\frac{m}{2}&\equiv& \frac{m}{2}\cdot \delta_{-1,\varepsilon}+k\mod m\Z
\qquad\textup{for } m\in F_2.
\end{eqnarray*}
If $\varepsilon=-1$ then any $m\in F_1$ is even so $l(m,2)\geq 1$.
In any case, whether $\varepsilon=-1$ or $\varepsilon=1$,
\begin{eqnarray}\label{4.10}
0\equiv 2^{l(m,2)-1}\cdot \delta_{-1,\varepsilon}+k\mod 2^{l(m,2)}\Z
\qquad\textup{for }m\in F_1.
\end{eqnarray}
Observe that $F_1$ contains elements $m_1\in E_1$ and $m_2\in E_2$
with $l(m_1,2)>l(m_2,2)$ as there is a highest 2-edge from $E_1$ to $E_2$.
This and \eqref{4.10} show $\varepsilon=1$. Now
\begin{eqnarray*}
0&\equiv& k\mod 2^{l(m,2)}\Z
\qquad\textup{for }m\in F_1\\
2^{l(m,2)-1}&\equiv& k\mod 2^{l(m,2)}\Z
\qquad\textup{for }m\in F_2
\end{eqnarray*}
follows. With $l(F_1,2)=l(M,2)$, the first congruence says
$2^{l(M,2)}|k$, the second congruence contradicts this, a contradiction.
Therefore $\varepsilon\in\{\pm 1\}$ and $k\in\Z$ with \eqref{4.9}
do not exist. One concludes as in the 1st case. $(S_2)$ is necessary
in theorem \ref{t1.2} in the case when $\GG(M)$ is connected.

\section{Sufficiency of the conditions in the main result
in the connected case}\label{c5}
\setcounter{equation}{0}

\noindent
The aim of this section is to show that the conditions in
case (I) in theorem \ref{t1.2} are sufficient for \eqref{1.6}
if $\GG(M)$ is connected.

Let $(H_M,h_M,S)$ be a triple as in the introduction and let
$M$ be the set of orders of $h_M$. Suppose that $\GG(M)$ is
connected and satisfies $(S_2)$ and $(T_p)$ for any 
prime number $p\geq 3$.
Let $c(x)\in\Z[x]$ be a polynomial with $|c(\lambda)|=1$
for any eigenvalue $\lambda$ of $h_M$. We want to show
that $\varepsilon\in\{\pm 1\}$ and $k\in\Z_{\geq 0}$ with
\begin{eqnarray}\label{5.1}
c(x)\equiv \varepsilon\cdot x^k \mod (\prod_{m\in M}\Phi_m).
\end{eqnarray}
exist. With lemma \ref{t4.1} this implies \eqref{1.6}.

First, a sign $\alpha(m)\in\{\pm 1\}$ and a number
$a(m)\in\{0,1,...,m-1\}$ are associated 
to any $m\in M$ be requiring
\begin{eqnarray}\label{5.2}
c(e(\frac{1}{m}))&=& \alpha(m)\cdot e(\frac{a(m)}{m})\\
\textup{and additionally }\quad \alpha(m)&=& 1\quad \textup{ if }m
\textup{ is even.}\label{5.3}
\end{eqnarray}
They exist and are unique by \eqref{3.4}.

Now we have to apply theorem \ref{t3.4} in order to link
the pairs $(\alpha(m),a(m))$ for different $m$ and for
varying prime numbers $p$.
This will prepare the choice of $\varepsilon\in\{\pm 1\}$ 
and $k\in\Z$ such that \eqref{5.1} holds.
We consider the same cases as in theorem \ref{t3.4}.
Lemma \ref{t5.1} is a straightforward application of it.

\begin{lemma}\label{t5.1}
(a) Let $p\geq 3$ be a prime number. 
Suppose that a $p$-edge goes from $m_1\in M$ to $m_2\in M$.
Then
\begin{eqnarray}\label{5.4}
\alpha(m_1)&=&\alpha(m_2),\\
a(m_1)&\equiv& a(m_2)\mod m_2\Z.\label{5.5}
\end{eqnarray}

(b) Suppose that a 2-edge goes from $m_1\in M$ to $m_2\in M$.
Then $\alpha(m_1)=1$ (by definition) and 
\begin{eqnarray}\label{5.6}
a(m_1)&\equiv& \frac{m_2}{2}\cdot \delta_{-1,\beta(m_1,m_2)}+a(m_2)\mod m_2\Z\\
&&\textup{for some }\beta(m_1,m_2)\in\{\pm 1\}
\textup{ if }m_2\textup{ is even}\nonumber,\\
a(m_1)&\equiv& a(m_2)\mod m_2\Z\quad 
\textup{ if }m_2\textup{ is odd}.\label{5.7}
\end{eqnarray}

(c) Let $m_1,m_2,m_3\in M$ be such that a 2-edge goes from
$m_1$ to $m_2$ and a 2-edge goes from $m_2$ to $m_3$.
Then $\alpha(m_1)=\alpha(m_2)=1$ (by definition) and 
\begin{eqnarray}
\beta(m_1,m_2)&=& \beta(m_1,m_3)=\beta(m_2,m_3)\quad
\textup{if }m_3\textup{ is even,}\label{5.8}\\
\beta(m_1,m_2)&=& \alpha(m_3)\qquad\textup{if }m_3\textup{ is odd.}
\label{5.9}
\end{eqnarray}
\end{lemma}

{\bf Proof:}
(a) Define 
$$c_2(x):=\alpha(m_1)\cdot x^{m_1-a(m_1)}\cdot c(x)\in\Z[x].$$
Then
\begin{eqnarray*}
c_2(e(\frac{1}{m_1}))&=&1,\\
c_2(e(\frac{1}{m_2}))&=&\alpha(m_1)\alpha(m_2)\cdot 
e(\frac{-a(m_1)+a(m_2)}{m_2})\\
&=& 1 \qquad (\textup{by theorem \ref{3.4} (a)}).
\end{eqnarray*}
If $m_1$ and $m_2$ are even then $\alpha(m_1)=\alpha(m_2)=1$
by definition.
If $m_1$ and $m_2$ are odd then $\alpha(m_1)\alpha(m_2)=1$ 
because $-1\notin\{e(\frac{k}{m_2})\, |\, k\in\Z\}$.
In any case \eqref{5.4} and \eqref{5.5} hold.

\medskip
(b) $m_1$ is even, thus $\alpha(m_1)=1$.
Define 
$$c_3(x):= x^{m_1-a(m_1)}\cdot c(x)\in\Z[x].$$
Then
\begin{eqnarray*}
c_3(e(\frac{1}{m_1}))&=&1,\\
c_3(e(\frac{1}{m_2}))&=&\alpha(m_2)\cdot 
e(\frac{-a(m_1)+a(m_2)}{m_2})\\
&=& \beta\quad\textup{ for some }\beta\in\{\pm 1\} 
\qquad (\textup{by theorem \ref{t3.4} (b)}).
\end{eqnarray*}
If $m_2$ is even, then $\alpha(m_2)=1$ and 
\begin{eqnarray*}
a(m_1)\equiv \frac{m_2}{2}\cdot\delta_{-1,\beta}+a(m_2)\mod m_2\Z.
\end{eqnarray*} 
If $m_2$ is odd, then $\alpha(m_2)=\beta$ because 
$-1\notin\{e(\frac{k}{m_2})\, |\, k\in\Z\}$, and then
\begin{eqnarray*}
a(m_1)\equiv a(m_2)\mod m_2\Z.
\end{eqnarray*}

\medskip
(c) $m_1$ and $m_2$ are even, thus $\alpha(m_1)=\alpha(m_2)=1$.
Define
$$c_4(x):=\beta(m_1,m_2)\cdot 
x^{(7+\beta(m_1,m_2))\cdot m_1/4-a(m_1)}\cdot c(x).$$
Then 
\begin{eqnarray*}
c_4(e(\frac{1}{m_1}))&=&1,\\
c_4(e(\frac{1}{m_2}))&=& \beta(m_1,m_2)\cdot 
e(\frac{-a(m_1)+a(m_2)}{m_2})=1 \qquad (\textup{by \ref{5.6}}),\\
c_4(e(\frac{1}{m_3}))&=&\beta(m_1,m_2)\alpha(m_3)\cdot 
e(\frac{-a(m_1)+a(m_3)}{m_3})\\
&=& 1 \qquad (\textup{by theorem \ref{t3.4} (c)}).
\end{eqnarray*}
If $m_3$ is even then $\alpha(m_3)=1$ and \eqref{5.6}
for $m_1$ and $m_3$ gives $\beta(m_1,m_2)=\beta(m_1,m_3)$.
If $m_3$ is odd then \eqref{5.7} for $m_1$ and $m_3$ 
gives $\beta(m_1,m_2)=\alpha(m_3)$.

As $m_3|\frac{m_2}{2}$, \eqref{5.6} also says
$a(m_1)\equiv a(m_2)\mod m_3\Z$. This shows 
$\beta(m_1,m_3)=\beta(m_2,m_3)$ if $m_3$ is even.
\hfill$\Box$

\begin{corollary}\label{t5.2}
Let $p$ and $q$ be prime numbers with $q\geq 3$
(here $p=q$ as well as $p\neq q$ are possible). Suppose that a 
$p$-edge goes from $m_1\in M$ to $m_2\in M$. Then
\begin{eqnarray}\label{5.10}
a(m_1)&\equiv& a(m_2)\mod q^{l(m_2,q)}\Z,\\
a(m_1)&\equiv& a(m_2)\mod 2^{l(m_2,2)}\Z\quad\textup{if }p\geq 3,\label{5.11}\\
a(m_1)&\equiv& a(m_2)\mod 2^{l(m_2,2)}\Z\quad\textup{if }p=2
\textup{ and }m_2\textup{ is odd},\label{5.12}\\
a(m_1)&\equiv& 2^{l(m_2,2)-1}\cdot\delta_{-1,\beta(m_1,m_2)}+ a(m_2)
\mod 2^{l(m_2,2)}\Z\label{5.13}\\
&&\textup{if }p=2 \textup{ and }m_2\textup{ is even}.\nonumber
\end{eqnarray}
\end{corollary}

{\bf Proof:}
If $p\geq 3$ \eqref{5.10} and \eqref{5.11} follow from \eqref{5.5}.
If $p=2$ \eqref{5.10} follows from \eqref{5.6} and \eqref{5.7}.
\eqref{5.12} follows from \eqref{5.7}.
\eqref{5.13} follows from \eqref{5.6}.
\hfill $\Box$ 

\bigskip

By hypothesis, $\GG(M)$ is connected and satisfies $(S_2)$ and 
$(T_p)$ for any prime number $p\geq 3$.
Therefore $(M,E(M)-\{\textup{highest 2-edges}\})$ has either
1 or 2 components. $(S_2)$ and remark \ref{t1.3} (iv) say 
about the two cases the following.

\medskip
{\bf Case (1),} there is only 1 component $M$: It is a single 2-plane.
Then choose $m_1^{(2)}\in M$ arbitrary.

\medskip
{\bf Case (2),} there are 2 components: 
One of them is the unique highest 2-plane $E_1$,
and the other component $E_2$ satisfies $l(E_2,2)<l(E_1,2)$.
Furthermore, there is a highest 2-edge from a vertex
$m_1^{(2)}\in E_1$ to a vertex $m_2^{(2)}\in E_2$.
Observe also that $l(m_1^{(2)},2)=l(E_1,2)=l(M,2)$ as $l(m,2)$ is constant 
for all vertices $m$ within one 2-plane.

\medskip
In both cases, choose for any prime number $q\geq 3$ a vertex
$m_1^{(q)}$ in the unique highest $q$-plane. 
Then $l(m_1^{(q)},q)=l(M,q)$. 
Now we define candidates $\varepsilon$ and $k$ which shall
satisfy \eqref{5.1}. Define
\begin{eqnarray}\label{5.14}
\varepsilon&:=& \left\{\begin{array}{lll}
\alpha(m_1^{(2)}) &\in\{\pm 1\} & \textup{ in case (1),}\\
\beta(m_1^{(2)},m_2^{(2)})&\in\{\pm 1\} & 
\textup{ in case (2) if }m_2^{(2)}\textup{ is even},\\
\alpha(m_2^{(2)})&\in\{\pm 1\} & 
\textup{ in case (2) if }m_2^{(2)}\textup{ is odd},
\end{array}\right. \\
k&\in&\Z_{\geq 0}\quad \textup{ such that }\nonumber \\
k&\equiv& a(m_1^{(q)})\mod q^{l(M,q)}
\quad\textup{for any prime number }q\geq 3,\label{5.15}\\
k&\equiv& \left\{\begin{array}{ll}
a(m_1^{(2)}) \mod 2^{l(M,2)} &\textup{in case (1),}\\
2^{l(M,2)-1}\cdot\delta_{-1,\varepsilon}+a(m_1^{(2)})  
\mod 2^{l(M,2)} &\textup{in case (2),}
\end{array}\right.  \label{5.16}
\end{eqnarray}
Here observe that for any prime number $p$ 
$l(m_1^{(p)},p)=l(M,p)$ because $m_1^{(p)}$ is in the
unique highest $p$-plane.
$k\in\Z$ can be chosen as in \eqref{5.15} and \eqref{5.16} 
because of the chinese remainder theorem.
We want to show that these $\varepsilon$ and $k$ satisfy \eqref{5.1}.

\medskip
{\bf Case (1):} Then $M$ is a single 2-plane, there are no 2-edges,
and $l(m,2)=l(M,2)$ for any $m\in M$. 
As $\GG(M)$ is connected, \eqref{5.16} and \eqref{5.11} imply
\begin{eqnarray}\label{5.17}
k\equiv a(m)\mod 2^{l(m,2)}\qquad\textup{for any }m\in M.
\end{eqnarray}
As $\GG(M)$ is connected, \eqref{5.14} and \eqref{5.4} imply
\begin{eqnarray}\label{5.18}
\varepsilon = \alpha(m)\qquad\textup{for any }m\in M.
\end{eqnarray}
Let $q\geq 3$ be a prime number.
As $\GG(M)$ is connected and satisfies $(T_q)$
(compare the remarks \ref{1.3} (ii) and (iii)),
\eqref{5.15} and \eqref{5.10} imply
\begin{eqnarray}\label{5.19}
k\equiv a(m)\mod q^{l(m,q)}\qquad\textup{for any }m\in M.
\end{eqnarray}
Together \eqref{5.17} and \eqref{5.19} give
\begin{eqnarray}\label{5.20}
k\equiv a(m)\mod m\qquad\textup{for any }m\in M.
\end{eqnarray}
Together \eqref{5.18} and \eqref{5.20} and \eqref{5.2} say
\begin{eqnarray}\label{5.21}
c(e(\frac{1}{m})) = \varepsilon\cdot e(\frac{k}{m}).
\qquad\textup{for any }m\in M.
\end{eqnarray}
This implies \eqref{5.1}.

\medskip
{\bf Case (2):} 
As $\GG(M)$ is connected and satisfies $(T_q)$  
for any prime number $q\geq 3$ 
(compare the remarks \ref{1.3} (ii) and (iii)),
\eqref{5.15} and \eqref{5.10} imply 
for any prime number $q\geq 3$
\begin{eqnarray}\label{5.22}
k\equiv a(m)\mod q^{l(m,q)}\qquad\textup{for any }m\in M.
\end{eqnarray}
Below we will show inductively 
\begin{eqnarray}\label{5.23}
k\equiv 2^{l(m,2)-1}\cdot \delta_{-1,\varepsilon}+ a(m)
\mod 2^{l(m,2)}\quad\textup{for even }m\in M,\\
\left. \begin{array}{lll}
k&\equiv& a(m)\mod 2^{l(m,2)}\\
\varepsilon&=&\alpha(m) \end{array}\right\}
\quad \textup{for odd }m\in M. \label{5.24}
\end{eqnarray}
Together \eqref{5.22}, \eqref{5.23} and \eqref{5.24} give
\begin{eqnarray}\label{5.25}
k\equiv \frac{m}{2}\cdot \delta_{-1,\varepsilon}+ a(m)
\mod m\quad\textup{for even }m\in M,\\
\left. \begin{array}{lll}
k&\equiv& a(m)\mod m \\
\varepsilon&=&\alpha(m) \end{array}\right\}
\quad \textup{for odd }m\in M .\label{5.26}
\end{eqnarray}
Together \eqref{5.25}, \eqref{5.26} and \eqref{5.2} say
\begin{eqnarray}\label{5.27}
c(e(\frac{1}{m})) = \varepsilon\cdot e(\frac{k}{m}).
\qquad\textup{for any }m\in M.
\end{eqnarray}
This implies \eqref{5.1}.

Therefore it rests for the proof in case (2) to show
\eqref{5.23} and \eqref{5.24}. 
The proof of \eqref{5.23} and \eqref{5.24} will consist of two
inductions.
The first induction will show the following slightly weaker statements:
\begin{eqnarray}\label{5.28}
k&\equiv& 2^{l(m,2)-1}\cdot \delta_{-1,\gamma(m)}+ a(m)
\mod 2^{l(m,2)}
\textup{ for even }m,\\
k&\equiv& a(m)\mod 2^{l(m,2)}\qquad\textup{for odd }m\in M,
\label{5.29}
\end{eqnarray}
with a unique $\gamma(m)\in\{\pm 1\}$ for even $m\in M$.
The second induction will show $\gamma(m)=\varepsilon$
for even $m\in M$ and $\alpha(m)=\varepsilon$ for odd $m\in M$.
This and \eqref{5.28} and \eqref{5.29} give \eqref{5.23} and \eqref{5.24}.

\medskip
{\bf The first induction:}
$(S_2)$ and $\GG(M)$ connected imply $(T_2)$, see remark \ref{t1.3} (iv).
Therefore starting at $m_1^{(2)}$, one can reach any $m\in M$ 
going through a chain of edges, in correct direction through
2-edges and in any direction through $p$-edges for $p\geq 3$.
$m_1^{(2)}$ satisfies \eqref{5.28} with $\gamma(m_1^{(2)})=\varepsilon$
by \eqref{5.16}. 

If a $p$-edge for some $p\geq 3$ goes from $m_1$ to $m_2$ and one of
them satisfies \eqref{5.28} or \eqref{5.29}, then the other satisfies
\eqref{5.28} or \eqref{5.29} too, and $\gamma(m_1)=\gamma(m_2)$.
This follows from \eqref{5.11} and $l(m_1,2)=l(m_2,2)$.

If a 2-edge goes from $m_1$ to $m_2$ and $m_1$ satisfies \eqref{5.28},
then $m_2$ satisfies \eqref{5.29} if $m_2$ is odd, because of \eqref{5.12}.
It satisfies \eqref{5.28} with $\gamma(m_2)=\beta(m_1,m_2)$
if $m_2$ is even, because of \eqref{5.13}.

This finishes the inductive proof of \eqref{5.28} and \eqref{5.29},
and it gives some additional information on $\gamma(m)$:
All elements $m$ in one 2-plane have the same value $\gamma(m)$
if they are even. If $m_2$ is even and is 
at the end of a 2-edge which starts at
$m_1$ then $\gamma(m_2)=\beta(m_1,m_2)$.

\medskip
{\bf The second induction:} We saw already $\gamma(m_1^{(2)})=\varepsilon$.
Therefore $\gamma(m)=\varepsilon$ for all $m$ in the 2-plane $E_1$.
\eqref{5.14} gives 
$\varepsilon=\beta(m_1^{(2)},m_2^{(2)})=\gamma(m_2^{(2)})$
if $m_2^{(2)}$ is even and
$\varepsilon=\alpha(m_2^{(2)})$ if $m_2^{(2)}$ is odd.
$(S_2)$ says that all 2-planes in $E_2$ are connected by 2-edges which
are not highest 2-edges. If there is a 2-edge from $m_1$ to $m_2$
and a 2-edge from $m_2$ to $m_3$ then \eqref{5.8} and \eqref{5.9}
show:
\begin{eqnarray*}
\textup{If }m_3\textup{ is even:}
&& \gamma(m_2)=\varepsilon\iff \gamma(m_3)=\varepsilon.\\
\textup{If }m_3\textup{ is odd:}
&& \gamma(m_2)=\varepsilon\iff \alpha(m_3)=\varepsilon.
\end{eqnarray*}
Therefore $\gamma(m)=\varepsilon$ for all even $m\in M$
and $\alpha(m)=\varepsilon$ for all odd $m\in M$.
This finishes the second induction and the proof of \eqref{5.23}
and \eqref{5.24} and the discussion of case (2).
\hfill$\Box$

\section{The proof of the main result in the disconnected case}\label{c6}
\setcounter{equation}{0}

\noindent
The aim of this section is to prove theorem \ref{t1.2}
in the case when $\GG(M)$ is not connected.
But before, we consider a more general situation and state a lemma.

For any unitary polynomial $f\in \C[x]$ of degree 
$\deg f\geq 1$, let
\begin{eqnarray}\label{6.1}
H_f&:=& \Z[x]/(f)\cong\Z^{\deg f}\\
h_f&:=& \textup{multiplication by }x: H_f\to H_f.\nonumber
\end{eqnarray}
Then $(H_f,h_f)$ is a $\Z$-lattice $H_f$ of rank $\deg f$ with a cyclic
automorphism $h_f$, i.e. 
\begin{eqnarray}\label{6.2}
H_f =\bigoplus_{i=0}^{n-1}\Z\cdot h_f^i(e_1)
\end{eqnarray}
for some $e_1\in H_f$. In fact, here one can choose $e_1=1$.

Part (b) of the following lemma is a kind of 
chinese remainder theorem for such pairs.

\begin{lemma}\label{t6.1}
Let $f,g\in\Z[x]$ be unitary polynomials of degrees $\geq 1$. 

(a) $H_{fg}$ contains a unique primitive sublattice which is
$h_{fg}$-invariant and such that the characteristic polynomial
of $h_{fg}$ on it is $f$. It is $(g)/(fg)\subset H_{fg}$,
and $((g)/(fg),h_{fg})\cong (H_f,h_f)$.

(b)
\begin{eqnarray}\label{6.3}
(H_{fg},h_{fg})\cong (H_f,h_f)\times (H_g,h_g)
\iff |R(f,g)|=1.
\end{eqnarray}
If this holds then 
\begin{eqnarray}
\Aut(H_{fg},h_{fg})\cong \Aut(H_f,h_f)\times  \Aut(H_g,h_g).\label{6.4}
\end{eqnarray}
\end{lemma}

{\bf Proof:} 
(a) Over $\Q$ instead of $\Z$, the equality
\begin{eqnarray*}
(\ker f(h_{fg}):H_{fg}\otimes_\Z\Q\to H_{fg}\otimes_\Z\Q)
=(g)/(fg)\otimes_\Z\Q
\end{eqnarray*}
is obvious. Of course, there is a unique primitive sublattice
$U$ of $H_{fg}$ with $U\otimes_\Z\Q=(g)/(fg)\otimes_\Z\Q$.
It is $U=((g)/(fg)\otimes_\Z\Q)\cap H_{fg}$. 

Consider the isomorphism of $\Z$-lattices
\begin{eqnarray}\label{6.5}
\Phi: \Z[x]_{\leq \deg(fg)-1}\to H_{fg},\quad a(x)\mapsto 
a(x)\mod (fg).
\end{eqnarray}
Then $\Phi^{-1}((g)/(fg))=\Z[x]_{\leq deg f-1}\cdot g$ is a
primitive sublattice of $\Z[x]_{\leq \deg (fg)-1}$, as
\begin{eqnarray*}
\Z[x]_{\leq \deg(fg)-1}= \Z[x]_{\leq \deg f-1}\cdot g \oplus
\Z[x]_{\leq \deg g-1}.
\end{eqnarray*}
Therefore $(g)/(fg)$ is a primitive sublattice of $H_{fg}$.
The monodromy $h_{fg}$ on it is cyclic with generator $g\mod (fg)$.
This shows the isomorphism $((g)/(fg),h_{fg})\cong (H_f,h_f)$.

\medskip
(b) In $(H_f,h_f)\times (H_g,h_g)=(H_f\times H_g,h_f\times h_g)$, 
$(H_f,h_f)$ and $(H_g,h_g)$ are primitive $h_f\times h_g$-invariant
sublattices such that the characteristic polynomial of
$h_f\times h_g$ on them is $f$ respectively $g$.
Together they generate the full $\Z$-lattice $H_f\times H_g$.

In $(H_{fg},h_{fg})$, the sum of the primitive 
sublattices $(g)/(fg)$ and $(f)/(fg)$ is $(f,g)/(fg)$.
It is a sublattice of $H_{fg}$ of full rank $\deg (fg)$
if and only if $R(f,g)\neq 0$ by lemma \ref{t2.3} (b).
Then it is a direct sum of the sublattices $(g)/(fg)$ and
$(f)/(fg)$, and then it is of index $|R(f,g)|$ in $H_{fg}$ 
by \eqref{2.10}.

Therefore $(H_{fg},h_{fg})\cong (H_f\times H_g,h_f\times h_g)$
if and only if $|R(f,g)|=1$.

\eqref{6.4} is an immediate consequence of \eqref{6.3}.
\hfill$\Box$ 

\bigskip

Now we return to the situation in section \ref{c1}.
The following elementary observations will be useful.

\begin{remarks}\label{t6.2}
Let $M\subset \Z_{\geq 1}$ be a finite set of positive integers.
Let $(H_M,h_M,S)$ be a triple as in section \ref{c1} such 
that $M$ is the set of orders of the eigenvalues of $h_M$.
\begin{eqnarray}\label{6.6}
-\id\in\{h_M^k\, |\, k\in\Z\}\iff 
\exists\ l\in\Z_{\geq 1}\textup{ mit }
\forall\ m\in M\ l(m,2)=l.
\end{eqnarray}
\begin{eqnarray}
|\{\pm h_M^k\, |\, k\in\Z\}|=\left\{\begin{array}{ll}
\lcm(M)&\textup{if }\exists\ l\in\Z_{\geq 1}\textup{ with }
\forall\ m\in M\\
& l(m,2)=l,\\
2\lcm(M)&\textup{else.} \end{array}\right. \label{6.7}
\end{eqnarray}
\end{remarks}

{\bf Necessity of the conditions in case (II) in theorem \ref{t1.2}.}
Let $M\subset\Z_{\geq 1}$ be a finite set of positive integers,
and let $(H_M,h_M,S)$ be a triple as in section \ref{c1} such that
$M$ is the set of orders of the eigenvalues of $h_M$.
Suppose that $\GG(M)$ is not connected and that
$\Aut(H_M,h_M,S)=\{\pm h_M^k\, |\, k\in\Z\}$.
We will show that all conditions in case (II) in theorem \ref{t1.2} hold.

Let $M_1,...,M_r$ with $r\geq 2$ be the components of $\GG(M)$.
Define 
$$f_j:=\prod_{m\in M_j}\Phi_m\quad\textup{for }j=1,...,r.$$
As there are no edges between different components, \eqref{3.6} gives
for the resultants
\begin{eqnarray*}
|(R(f_i,f_j)|&=&1\quad\textup{ for }i\neq j,\\
|R(f_1,\prod_{j\geq 2}f_j)|&=&1,\ |R(f_2,\prod_{j\geq 3}f_j)|=1,..., 
\ |R(f_{r-1},f_r)|=1.
\end{eqnarray*}
One applies $r-1$ times lemma \ref{t6.1} and obtains
\begin{eqnarray*}
(H_M,h_M)\cong (H_{M_1},h_{M_1})\times ... \times (H_{M_r},h_{M_r}).
\end{eqnarray*}
As all eigenspaces are one-dimensional and $S$ is $h_M$-invariant,
it is clear that $S$ and the automorphisms decompose accordingly,
\begin{eqnarray*}
(H_M,h_M,S)&\cong& (H_{M_1},h_{M_1},S_1)\times ... 
\times (H_{M_r},h_{M_r},S_r),\\
\Aut(H_M,h_M,S)&\cong& \Aut(H_{M_1},h_{M_1},S_1)\times ... 
\times \Aut(H_{M_r},h_{M_r},S_r).
\end{eqnarray*}
Recall
\begin{eqnarray*}
\Aut(H_M,h_M,S)&=&\{\pm h_M^k\, |\, k\in\Z\}\quad\textup{ by hypothesis},\\
\Aut(H_{M_j},h_{M_j},S_j)&\supset& \{\pm h_{M_j}^k\, |\, k\in\Z\},\\
|\{\pm h_{M_j}^k\, |\, k\in\Z\}|
&=& \left\{\begin{array}{ll}\lcm(M_j)& \textup{if }(*j)\textup{ holds,}\\
2\lcm(M_j)& \textup{ else,}\end{array}\right. \\
(*j)&:&  -\id\in \{h_{M_j}^k\, |\, k\in\Z\}\\
&\iff& \exists\ l_j\in\Z_{\geq 1}\textup{ with }
\forall\ m\in M_j\ l(m,2)=l_j,\\
&\iff& M_j\textup{ is a 2-plane with }l(M_j,2)\geq 1.
\end{eqnarray*}
The last equivalence holds because $M_j$ is a connected subgraph
of $\GG(M)$.

\medskip
{\bf First case,} $\exists\ l\in\Z_{\geq 1}\textup{ with }
\forall\ m\in M\ l(m,2)=l$: Then all $M_j$ are 2-planes with
$l(M_j,2)=l\geq 1$, so $(*j)$ holds. Thus
\begin{eqnarray*}
|\Aut(H_M,h_M,S)|&=&\lcm(M)=\lcm(\lcm(M_1),...,\lcm(M_r)),\\
|\Aut(H_{M_j},h_{M_j},S_j)|&\geq& |\{\pm h_{M_j}^k\, |\, k\in\Z\}|
=\lcm(M_j),\\
\lcm(\lcm(M_1),...,\lcm(M_r))&\geq& \prod_{j=1}^r\lcm(M_j).
\end{eqnarray*}
But $2|\lcm(M_j)$ for all $j$. This is a contradiction. 
The first case is impossible.

\medskip
{\bf Second case,} $\not\exists\ l\in\Z_{\geq 1}\textup{ with }
\forall\ m\in M\ l(m,2)=l$: Then
\begin{eqnarray*}
|\Aut(H_M,h_M,S)|&=&2\lcm(M)\\
&=&2\lcm(\lcm(M_1),...,\lcm(M_r)),\\
|\Aut(H_{M_j},h_{M_j},S_j)|&\geq& 
\left\{\begin{array}{ll}\lcm(M_j)& \textup{if }(*j)\textup{ holds,}\\
2\lcm(M_j)& \textup{ else,}\end{array}\right. \\
2\lcm(\lcm(M_1),...,\lcm(M_r))&\geq& 
\prod_{j=1}^r\left\{
\begin{array}{ll}\lcm(M_j)&\textup{if }(*j)\textup{ holds}\\ 
2\lcm(M_j)&\textup{else}\end{array}\right\}
\end{eqnarray*}
Therefore $|\{j\, |\, (*j)\textup{ holds}\}|=r-1$ (case (i)) 
or $=r$ (case (ii)).
In the case (i) $\gcd(\lcm(M_i),\lcm(M_j))=1$ for all $i\neq j$,
which implies $r=2$. We can suppose that $(*1)$ holds.
Then $M_1$ is a 2-plane with $l(M_1,2)\geq 1$. As $l(M_2,2)=0$,
also $M_2$ is a 2-plane.

In case (ii)
$\gcd(\lcm(M_i),\lcm(M_j))=1$ for all pairs $(i,j)$ with $i<j$
except for one pair which may satisfy $\gcd(\lcm(M_i),\lcm(M_j))=2$.
This also implies $r=2$. $M_1$ and $M_2$ are 2-planes with 
$l(M_i,2)\geq 1$ and $\gcd(\lcm(M_1),\lcm(M_2))=2$.
As we are in the second case, $l(M_1,2)\neq l(M_2,2)$.
We can suppose $l(M_1,2)> l(M_2,2)=1$. 

In case (i) as well as in case (ii), the inequalities above
are equalities, and thus 
$$|\Aut(H_M,h_M,S)|=|\{\pm h_M^k\, |\, k\in\Z\}|$$
and thus $\Aut(H_M,h_M,S)=\{\pm h_M^k\, |\, k\in\Z\}$. 
Case (I) in theorem \ref{t1.2} says that the graphs 
$\GG(M_1)$ and $\GG(M_2)$ satisfy $(T_p)$ for any prime number $p\geq 3$.
This completes the proof of the necessity of the conditions in case (II)
in theorem \ref{t1.2}.

\bigskip
{\bf Sufficiency of the conditions in case (II) in theorem \ref{t1.2}.}
Let $M\subset\Z_{\geq 1}$ be a finite set of positive integers,
and let $(H_M,h_M,S)$ be a triple as in section \ref{c1} such that
$M$ is the set of orders of the eigenvalues of $h_M$.
Suppose that $\GG(M)$ has two components $M_1$ and $M_2$
whose graphs are 2-planes and satisfy $(T_p)$ for any prime number
$p\geq 3$. Suppose also  \eqref{1.7} and \eqref{1.8}, which are
\begin{eqnarray*}
\gcd(\lcm(M_1),\lcm(M_2))&\in&\{1;2\},\\
l(M_1,2)&>& l(M_2,2)\in\{0;1\}.
\end{eqnarray*}
We want to show $\Aut(H_{M},h_{M},S)=\{\pm h_M^k\, |\, k\in \Z\}$.

As there are no edges between $M_1$ and $M_2$ in the graph $\GG(M)$,
\begin{eqnarray*}
|R(\prod_{m\in M_1}\Phi_m,\prod_{m\in M_2}\Phi_m)|=1.
\end{eqnarray*}
Lemma \ref{t6.1} applies and says
\begin{eqnarray*}
(H_M,h_M)\cong (H_{M_1},h_{M_1})\times (H_{M_2},h_{M_2}).
\end{eqnarray*}
As all eigenspaces are one-dimensional and $S$ is $h_M$-invariant,
$S$ and the automorphisms decompose accordingly,
\begin{eqnarray*}
(H_M,h_M,S)&\cong& (H_{M_1},h_{M_1},S_1)\times (H_{M_2},h_{M_2},S_2),\\
\Aut(H_M,h_M,S)&\cong& \Aut(H_{M_1},h_{M_1},S_1)\times \Aut(H_{M_2},h_{M_2},S_2).
\end{eqnarray*}
As $\GG(M_1)$ and $\GG(M_2)$ satisfy $(T_p)$ for any prime number
$p\geq 3$, case (I) of theorem \ref{t1.2} applies and gives
\begin{eqnarray*}
\Aut(H_{M_i},h_{M_i},S_i)=\{\pm h_{M_i}^k\, |\, k\in\Z\}\quad
\textup{for }i=1,2.
\end{eqnarray*}
$M_1$ and $M_2$ are 2-planes with 
$l(M_1,2)\geq 1$ and $l(M_2,2)\in\{0;1\}$, thus
\begin{eqnarray*}
|\Aut(H_{M_1},h_{M_1},S_1)|&=&\lcm(M_1),\\
|\Aut(H_{M_2},h_{M_2},S_2)|&=&\left\{\begin{array}{ll}
\lcm(M_1)&\textup{if }l(M_2,2)=1,\\
2\lcm(M_1)&\textup{if }l(M_2,2)=0.\end{array}\right. 
\end{eqnarray*}
Therefore 
\begin{eqnarray*}
&&|\Aut(H_{M},h_{M},S)|\\
&=&|\Aut(H_{M_1},h_{M_1},S_1)|\cdot|\Aut(H_{M_2},h_{M_2},S_2)|\\
&=&\left\{\begin{array}{ll}
\lcm(M_1)\cdot \lcm(M_2)&\textup{if }l(M_2,2)=1,\\
\lcm(M_1)\cdot 2\lcm(M_2)&\textup{if }l(M_2,2)=0,\end{array}\right. \\
&=& 2\cdot \lcm(M)
= |\{\pm h_M^k\, |\, k\in \Z\}|,
\end{eqnarray*}
thus
$$\Aut(H_{M},h_{M},S)=\{\pm h_M^k\, |\, k\in \Z\}.$$
\hfill$\Box$


\begin{thebibliography}{99}
\bibitem[Ap70]{Ap70} T. Apostol: \quad
     Resultants of cyclotomic polynomials.
     Proc. A.M.S. {\bf 24} (1970), 457--462.
\bibitem[GH16]{GH16} F. Gau{\ss}, C. Hertling: \quad 
    $\mu$-constant monodromy groups and Torelli results for 
    marked singularities, for the unimodal and some bimodal    
    singularities.   
    In: Singularities and Computer Algebra, Festschrift for Gert-Martin Greuel 
    on the Occasion of his 70th Birthday (W. Decker, G. Pfister, M. Schulze, eds.).
    Springer International Publishing 2017, 109-146. 
\bibitem[GH17]{GH17} F. Gau{\ss}, C. Hertling: \quad 
    $\mu$-constant monodromy groups and Torelli results for the 
    quadrangle singularities and the bimodal series.
    arXiv:1710.03507v1, 125 pages, 10.10.2017. 
\bibitem[He98]{He98} C. Hertling: \quad  
    Brieskorn lattices and Torelli type
    theorems for cubics in $\P^3$ and for Brieskorn-Pham 
    singularities
    with coprime exponents. In: Singularities, the Brieskorn 
    anniversary
    volume. Progress in Mathematics {\bf 162}. 
    Birkh\"auser Verlag, Basel-Boston-Berlin 1998, pp. 167--194.
\bibitem[He11]{He11} C. Hertling: \quad
    $\mu$-constant monodromy groups and marked singularities.
    Ann. Inst. Fourier, Grenoble {\bf 61.7} (2011), 2643--2680.
\bibitem[HZ18]{HZ18} C. Hertling, Ph. Zilke: \quad 
    Seven combinatorial problems around quasihomogeneous singularities.
    Manuskript, January 2018, 36 pages.
\bibitem[Or72]{Or72} P. Orlik: \quad 
    On the homology of weighted homogeneous polynomials. 
    In: Lecture Notes in Math. 298, Springer, Berlin, 1972.
\bibitem[vW71]{vW71} B.L. van der Waerden: \quad
    Algebra I, Springer, Berlin, Heidelberg, New York, 
    8. Auflage, 1971.
\bibitem[Wa82]{Wa82} L.C. Washington: \quad 
    Introduction to cyclotomic fields.
    Springer, New York, 1982.
\end{thebibliography}
\end{document}